\documentclass[11pt, leqno]{amsart}
\usepackage{amsthm, amsmath}
\usepackage{amssymb} 
\usepackage{amscd}
\usepackage[curve,matrix,arrow,frame,tips]{xy}
\usepackage{verbatim}
\usepackage{graphicx}
\newtheorem{thm}{Theorem}[section]
\newtheorem{prop}[thm]{Proposition}
\newtheorem{lemma}[thm]{Lemma}

\newtheorem{Remarknumb}[thm]{Remark}

\newtheorem{Remark}[thm]{Remark}

\newtheorem{conjecture}[thm]{Conjecture}
\newtheorem{cor}[thm]{Corollary}

\newcounter{ex}[section]

\newcommand{\cal}{\mathcal}

\newcommand{\ds}{\displaystyle}
\newcommand{\E}{{\mathcal E}}
\newcommand{\A}{{\mathcal A}}

\newcommand{\Bb}{{\mathcal B}}

\newcommand{\R}{{\bf R}}
\newcommand{\Q}{{\mathbb Q}}

\newcommand{\Ccc}{{\mathfrak C}}

\newcommand{\B}{{\mathcal B}}

\newcommand{\X}{{\mathcal X}}

\newcommand{\T}{{\mathcal T}}
\newcommand{\Ff}{{\Bbb F}}
\newcommand{\Gg}{{\mathcal G}}

\newcommand{\Gm}{{{\mathbb G}_m}}
\newcommand{\Z}{{\mathbb Z}}
\newcommand{\F}{{\mathcal F}}
\newcommand{\ti}{\tilde}
\newcommand{\Spec}{{\rm Spec\, }}
 \renewcommand{\O}{{\mathcal O}}

\newcommand{\e}{{\eta}}

\renewcommand{\r}{\theta}

\newcommand{\Mfr}{{\mathfrak M}}
\newcommand{\Sfr}{{\mathfrak S}}
\newcommand{\M}{{\mathcal M}}
\renewcommand{\L}{{\mathcal L}}

\newcommand{\Res}{{\rm Res}}

\newcommand{\RR}{{\mathcal R}}

\newcommand{\Cc}{{\mathcal C}}

\newcommand{\N}{{\mathcal N}}

\def\thfill{\null\nobreak\hfill}

\def\endproof{\thfill\vbox{\hrule
  \hbox{\vrule\hbox to 5pt{\vbox to 5pt{\vfil}\hfil}\vrule}\hrule}}

\newcommand{\Gal}{{\rm Gal}}

\renewcommand{\to}{\longrightarrow}

\newcommand{\xto}{\xrightarrow}

\begin{document}

\title[ ]{$\Phi$-modules and coefficient spaces}
\author[G. Pappas]{G. Pappas}
%\thanks{*Partially supported by  NSF Grant DMS05-01409.}
 
\address{Dept. of
Mathematics\\
Michigan State
University\\
E. Lansing\\
MI 48824-1027\\
USA}
\email{pappas@math.msu.edu}
\author[M. Rapoport]{M. Rapoport}
\address{Mathematisches Institut der Universit\"at Bonn,  
Beringstrasse 1\\ 53115 Bonn\\ Germany.}
\email{rapoport@math.uni-bonn.de}

\date{\today}
 
\maketitle

% \tableofcontents

 \medskip
 \section*{Introduction}

This paper is inspired by Kisin's article \cite{K1}, in which he studies deformations of Galois representations
of a local $p$-adic field which are defined by finite flat group schemes. The result of Kisin most relevant
to our paper is his construction of a kind of resolution of the formal deformation space of the given
Galois representation, by constructing a scheme which classifies all finite flat group schemes giving
rise to the deformed Galois representation. Our purpose here is to globalize Kisin's construction.

Let $K$ be a finite extension of $\Q_p$, with residue field $k$. Let $K_0$ be the maximal unramified
extension of $\Q_p$ contained in $K$. Then $K_0$ is the fraction field of the ring of Witt vectors
$W = W(k)$. Let $\pi$ be a uniformizer of $K$ and $E(u)\in W[u]$ the Eisenstein polynomial that
$\pi$ satisfies. Let $G_K = \Gal (\overline{K}/K)$ be the absolute Galois group of $K$.
Set $\Sfr$ for the ring of formal power series $W[[u]]$. Let $\phi : \Sfr\rightarrow\Sfr$ be such that $\phi\vert W$ is the Frobenius automorphism
and with $\phi (u) = u^p$.

Kisin's construction is based on the existence of a fully faithful exact functor from a suitable category of 
$\Sfr$-modules $\Mfr$ equipped with a $\phi$-linear endomorphism $\Phi$ to the category of finite flat
(commutative)
group schemes of $p$-power rank over $\Spec(\O_K)$. This in turn was inspired by work of Breuil \cite{Br} 
who gave a similar but more complicated description of such group schemes.
A variant of this functor also works with coefficients:
if $R$ is a $\Z_p$-algebra with finitely many elements, then there is a similar functor from a suitable
category of $\Sfr\otimes_{\Z_p} R$-modules $\Mfr$ with $\phi$-linear endomorphism $\Phi$ to the
category of finite flat group schemes of $p$-power rank with $R$-action over $\Spec\O_K$.

Let $K_\infty / K$ be the extension obtained by adjoining a compatible system of $p^n$-power roots of
$\pi$, and let $G_{K_\infty} = \Gal (\overline{K} / K_\infty)$ be its absolute Galois group. Let $\O_{\E}$
be the $p$-adic completion of $\Sfr [1/u]$, a complete discrete valuation ring, with uniformizer
$p$ and residue field $k((u))=k[[u]][1/u]$. Then there exists an equivalence of categories between the category of
finitely generated $\O_\E$-modules $M$ equipped with an isomorphism $\Phi : \phi^* (M)\rightarrow M$
and the category of continuous representations of $G_{K_\infty}$ in $\Z_p$-modules, and this is compatible with
the previous functor via the restriction functor from $G_K$-representations to $G_{K_\infty}$-representations.
Again there is also a variant for representations with values in a finite coefficient $\Z_p$-algebra $R$.

Our basic idea in this paper is to formally forget about Galois representations and finite flat group
schemes and simply consider the modules themselves, without any finiteness conditions on $R$.
More precisely, for any $\Z_p$-algebra $R$ set $R_W = W\otimes_{\Z_p} R$ and extend the
endomorphism $\phi$ of $W((u))=W[[u]][1/u]$ to $R_W ((u)) = W ((u))\hat{\otimes}_{\Z_p} R$ by the
identity on the second factor. We define the fpqc-stack $\Cc_{d, K}$ by giving its values on
$\Z_p$-algebras $R$ as the groupoid of pairs $(\Mfr, \Phi)$, where $\Mfr$ is a finitely generated
$R_W [[u]]$-module which is free of rank $d$ locally fpqc on $\Spec R$ and where
\begin{equation*}
\Phi : \phi^*  \Mfr [1/u]\xto{\sim}\Mfr [1/u]
\end{equation*}
is an isomorphism of $R_W ((u))$-modules such that $E (u) \Mfr\subset\Phi (\phi^* \Mfr)\subset\Mfr$.

We also introduce the fpqc-stack $\mathcal R_d$ with values in a $\Z_p$-algebra $R$ the groupoid of pairs
$(M, \Phi)$, where $M$ is a finitely generated $R_W ((u))$-module which is locally fpqc on $\Spec R$
free of rank $d$ as $R_W ((u))$-module, and where $\Phi :\phi^*  (M)\xto{\sim} M$.

There is an obvious morphism
\begin{equation*}
\r : \Cc_{d, K}\rightarrow\mathcal R_d
\end{equation*}
sending $(\Mfr, \Phi)$ to $(\Mfr [1/u], \Phi)$.

Our main results concern the algebraicity of the previous construction. Let $\widehat\Cc_{d, K}$ be
the $p$-adic completion of the stack $\Cc_{d, K}$, i.e. its restriction to $p$-nilpotent
$\Z_p$-algebras. Then
\begin{equation*}
\widehat\Cc_{d, K} = \varinjlim\nolimits_a \Cc_{d, K}\times_{\Spec \Z_p}\Spec \Z / p^a\Z\ .
\end{equation*}
Our main result shows that this presents $\widehat\Cc_{d, K}$ as an inductive $2$-limit of Artin
stacks of finite type over $\Spec\Z /p^a\Z$. Furthermore, the singularities of $\widehat\Cc_{d, K}$ are modeled by local models. 

\begin{thm}\label{thmintro}
(i) For each $a$, the stack $\Cc_{d, K}\times_{\Spec \Z_p}\Spec \Z / p^a\Z$ on $\Z / p^a\Z$-algebras
is representable by an Artin stack ${\Cc^a_{d, K}}$ of finite type over $\Spec\Z / p^a\Z$. The
inductive limit $\varinjlim_{a}{\Cc^a_{d, K}}$ is the formal $p$-adic completion $\widehat
\Cc_{d, K}$
of $\Cc_{d, K}$.

\smallskip

\noindent (ii) There is a ``local model" diagram
$$\xymatrix{
{} & \ar[dl]_\pi \widetilde{\widehat\Cc_{d, K}} \ar[dr]^\varphi & {}\\
\widehat\Cc_{d, K} & {} & \widehat M_{d, K}\ ,
 }$$
in which $\pi$ is a principal homogeneous space under the positive loop group $L^+ G$ of $G = \Res_{W / \Z_p} (GL_d)$ 
completed along its special fiber,  and in which $\varphi$ is formally
smooth. Here $M_{d, K}$ is the projective $\Z_p$-scheme parametrizing all $\O_K\otimes_{\Z_p}\O_S$-submodules
of $(\O_K\otimes_{\Z_p}\O_S)^d$ which are locally direct summands as $\O_S$-modules, 
and $\widehat M_{d, K}$ denotes  its completion along the special fiber.

\smallskip

\noindent (iii) For each $a$, set ${\mathcal R}^a_d={\mathcal R}_d\times_{\Spec \Z_p}\Spec \Z / p^a\Z$.  
The morphism $\r^a : {\Cc^a_{d, K}}\rightarrow {\mathcal R}^a_d$  given by reducing \ $\r$ modulo $p^a\Z$ is representable and 
proper; hence $\hat \r : \widehat\Cc_{d, K}\rightarrow \widehat {\mathcal R}_d$ is an inductive  limit  
of representable and proper morphisms. 
\end{thm}

The fibers of $\theta$ over a finite field $\Ff$  are interesting projective subvarieties of the 
affine Grassmannian of the group $G = \Res_{W / \Z_p} (GL_d)$, which we call {\it Kisin varieties}. More
precisely, we define variants of $\Cc_{d, K}$ and $M_{d, K}$ depending on a co-character $\mu$ of
$G$ and define Kisin varieties associated to $(G, A, \mu)$, where $A\in  G(\Ff((u)))$ defines the given
$\Ff$-valued point in $\mathcal R_d$. They are the analogues, for the kind of Frobenius involved here,  of the affine Deligne-Lusztig varieties appearing in the isocrystal context, cf., eg. \cite{ghkr}. The study of these varieties was begun by Kisin in \cite{K1}, in the case
$d = 2$ and $W = \Z_p$.  In a companion paper to ours, E. Hellmann extends Kisin's results (again for
$d = 2$ and $W = \Z_p$). In Hellmann's paper, one of the main tools is the Bruhat-Tits building of $GL_2$. We show here how the Bruhat-Tits building can be used in general to gain a qualitative overview of Kisin varieties.

The other extreme to the fiber over $\Ff$ of $\widehat \Cc_{d, K}$ is its fiber over $\Q_p$. Here we construct a kind of period map of stacks over the category of adic formal schemes  locally of finite type over ${\rm Spf} (\Z_p)$, 
$$
\Pi(\mathcal X): \widehat\Cc_{d, K}(\mathcal X)\rightarrow \mathfrak D_{d, K}(\mathcal X^{\rm rig})\ ,
$$
where $\mathfrak D_{d, K}$ is the stack over the category of rigid-analytic spaces over $\Q_p$ parametrizing filtered $\Phi$-modules (both the filtration and $\Phi$ vary!). To determine the image of the period map seems one of the major challenges in the theory. We conjecture that the image should be given somewhat analogously to Hartl's {\it admissible set} in \cite{Ha}. 

As is apparent from the above, the theory developed here bears many similarities to the theory of period spaces for $p$-divisible groups in \cite{R-Z}, but there are also substantial differences. The stack $\widehat \Cc_{d, K}$ is analogous to one of the period spaces of $p$-divisible groups in \cite{R-Z}, but unlike these it is adic over $\Spec(\Z_p)$ ($p$ generates an ideal of definition); the local model diagram looks formally just like the corresponding one in \cite{R-Z}; Kisin varieties are the analogues of affine Deligne-Lusztig varieties, and the stack ${\mathfrak D}_{d, K}$ plays a role similar to  the Grassmannian containing the period space of \cite{R-Z}. In \cite{R-Z}, the base scheme of the $p$-divisible groups is variable; here the base scheme $\Spec(\O_K)$ is constant, but the coefficients are variable.

We now explain the lay-out of the paper. In section 1 we explain by analogy on the classical theory of unit root crystals the spaces/stacks we encounter. In section 2 we prove our main technical result, which states that the stack $\Cc_d$, which associates to  $\Z_p$-algebras $R$ the groupoid of locally free $R_W[[u]]$-modules of rank $d$ with $\Phi$-module structure,  can be presented as an inductive limit of Artin stacks of finite type over $\Z_p$. In section 3 we fix the local field as above and prove the main theorem stated above. In section 4 we indicate the relation to the deformation spaces of Galois representations which is at the origin of Kisin's theory. In section 5 we construct and discuss the period morphism. In the final section 6 we discuss Kisin varieties and their analysis through Bruhat-Tits buildings. 

We thank X.\ Caruso, E.\ Hellmann, M.\ Kisin, R.\ Kottwitz, A.\ Lytchak, J.\ Stix and E.\ Viehmann for helpful discussions.

\bigskip
\medskip
 
\section{Motivation: Unit crystals and Galois representations} \label{1}

\subsection{Unit root crystals}  
Let $k$ be a finite field of characteristic $p>0$; for simplicity, we will assume that $k={\Ff}_p$. Let $X$ a variety over $k$. Denote by $\phi:
X\rightarrow X$ the Frobenius $\phi(a)=a^p$ for $a\in \O_X$. Suppose that $S$ is a $k$-scheme, and let $\phi_S=\phi\times id_S: X\times S\rightarrow X\times S$. Consider pairs $(\M, F)$ consisting  of a  locally free $\O_{X\times S}$-coherent sheaf $\M$ of rank $d$ on $X\times S$ and an isomorphism
$$
F: \phi^*_S\M\xrightarrow{\sim} \M\ .
$$
 As $S$ varies, these pairs form  an fpqc stack
$FM^{d, et}_X$ over $\Spec k$. In fact,  $FM^{d, et}_X$ is an Artin stack locally of finite type over $k$. Indeed, let ${ Fib}^d_{X/k}$ be the Artin stack locally of finite type over $k$, whose values in a $k$-scheme  $S$ is the groupoid of locally free $\O_{X\times S}$-modules of rank $d$, and let ${\widetilde{ Fib}}^d_{X/k}$ be the $GL_d$-torsor over ${ Fib}^d_{X/k}$, consisting of a locally free $\O_{X\times S}$-modules $\M$ of rank $d$ and a basis $\iota: \M \xrightarrow{\sim} \O^d_{X\times S}$. Then there is an action of $GL_d$ on the product ${\widetilde{ Fib}}^d_{X/k}\times GL_d$ via 
$$
g:\ \ (\M, \iota, A)\mapsto (\M, g^{-1}\cdot \iota, g^{-1}\cdot A\cdot \phi(g))\ .
$$
This presents $FM^{d, et}_X$ as a quotient of ${\widetilde{ Fib}}^d_{X/k}\times GL_d$ by $GL_d$, and hence $FM^{d, et}_X$ is an Artin stack locally of finite type, as claimed. 

\medskip

Suppose $S=\Spec(\Lambda)$ with $|\Lambda|<\infty$.
Then by Katz \cite{katz}, 4.1, cf.\ also  \cite{emerton-kisin},  there is a bijective correspondence between pairs $(\M, F)$ over $\Spec(\Lambda)$ and \'etale sheaves of $\Lambda$-modules on $X$. In the case $\Lambda={\Ff}_p$, this correspondence is obtained via push-out from the injection $GL_d({\Ff}_p)\rightarrow GL_d$ which induces an equivalence of categories between the category of $GL_d({\Ff}_p)$-torsors on $X$ (for the \'etale topology) and the category of $GL_d$-torsors $P$ on $X$ with an isomorphism $\phi^*(P)\xrightarrow{\sim} P$. 

\medskip

We want to think of $FM^{d, et}_X$ as a ``coefficient space" for $p$-torsion representations of $\pi_1(X, \bar\eta)$. However, it seems that global questions on these spaces (i.e., when $S$ is not a local  Artin ring)  have not been studied much in the literature. For instance, are there non-constant morphisms of projective $k$-schemes into $FM^{d, et}_X$? What is the dimension of $FM^{d, et}_X$? Etc. The only result we are aware of is Laszlo's construction \cite{La} of a projective curve $X$ of genus $2$ over the field with $2$ elements, a projective curve $S$ over a finite extension $k'$ of ${\Ff}_2$ and a locally free coherent $\O_{X\times S}$-module $\M$ of arbitrary rank with an isomorphism $F: (\phi_S^2)^*\M\xrightarrow{\sim} \M$.

\subsection{Variants} We mention here some variants of the above theory. Let $G$ be a reductive group over ${\Ff}_p$. Then we can consider the fpqc stack $FM^{G, et}_X$ of pairs $(P, F)$, where $P$ is a $G$-torsor on $X\times S$ and where $F: \phi_S^*(P)\xrightarrow{\sim}P$. For $G=GL_d$, we recover the stack considered above. 

We may also consider ``meromorphic Frobenius structures", as follows. Assuming $X$ to be irreducible, with generic point $\eta(X)$, consider the fpqc stack  $FM^d_X$ of pairs
$(\M, F)$ with $\M$ a locally free $\O_{X\times S}$-coherent sheaf of rank $d$ on $X\times S$ and 
$$
F: \phi^*_S\M\xrightarrow{ } \M
$$
a  homomorphism such that $F_{|\eta(X)\times S}$ is an isomorphism.

One may also control the degeneracy of the meromorphic Frobenius structure. For instance, let $X$ be a curve. Then we may consider triples $(\M, F, x)$ with $(\M, F)$ 
in $FM^d_X$ and $x: S\rightarrow X$ such that $Coker(F)$
is supported on the graph $\Gamma_x\subset X\times S$ 
and is annihilated by the power of the ideal sheaf $ I_{\Gamma_x}^e$ for some 
fixed $e\geq 1$. Denoting the corresponding stack by $FM^d_{X, e}$,
there is a morphism (the "pole morphism"), 
$$
p: FM^d_{X, e}\rightarrow X
$$

 Similarly we can obtain a construction that resembles 
shtuka, but the Frobenius is ``on the other factor". Namely, assume that 
$X$ is a curve as before. Consider $(\M, \M', F, F', x, y)$ with $\M$, $\M'$ 
locally free $\O_{X\times S}$-coherent sheaves of rank $d$ on $X\times S$ and homomorphisms
$$
\phi^*_S\M\xrightarrow{F} \M'\xleftarrow{F'} \M
$$
such that $Coker(F)$, resp. $Coker(F')$ 
is supported on the graph of $x: S\rightarrow X$, resp. 
$y: S\rightarrow X$. Again we can ask that $Coker(F)$, resp. $Coker(F')$ ,
satisfy some additional property.

Another variant is obtained by replacing the variety $X$ by the spectrum of the completed local ring at a closed point, or by its fraction field. 

\medskip

All these ``spaces"/stacks seem interesting geometric objects. 
\bigskip
\medskip

\section{Spaces of Kisin-Breuil modules} \label{2}

Fix a finite field $k$ of characteristic $p$ and denote by 
$\phi(a)=a^p$ the Frobenius automorphism of $k$. We will denote by $W=W(k)$ the ring of Witt vectors of $k$ and by $\phi : W\rightarrow W$ the unique lifting of Frobenius. 

Let $R$ be a commutative $\Z_p$-algebra and set $R_W=W\otimes_\Z R$. We extend $\phi$ in a $R$-linear way to $R_W$ and denote this extension also by $\phi$. We also denote by $\phi$ the 
endomorphism $\phi$ of $R_W((u))=W((u))\hat\otimes_{\Z} R$ 
 given by 
$$
\phi(\sum_i a_iu^i)=\sum_i \phi(a_i) u^{pi}\ .
$$

\subsection{}\label{2defofstacks}
We  define now various stacks of modules with Frobenius structure.

\medskip

Let us consider the stack $\Cc_{ d}$ such that $\Cc_d(R)$ 
is the groupoid of $R_W[[u]]$-$\Phi$-modules $(\Mfr, \Phi)$: These are by definition pairs
of a $R_W[[u]]$-module $\Mfr$ which is locally on $R$ (for the fpqc topology)
$R_W[[u]]$-free of rank $d$ and a 
$R_W((u))$-module isomorphism 
$$
\Phi : \phi^*{\frak M}[1/u]=R_W((u))\otimes_{\phi, R_W[[u]]}\Mfr\xrightarrow{\sim} {\Mfr}[1/u]=R_W((u))\otimes_{  R_W[[u]]}\Mfr\ . 
$$ 
It is easy to see that $\Cc_d$ is a stack for the fpqc topology.

\medskip

Next,  consider  the stack $\RR_d$ which is such that
$\RR_d(R)$ is the groupoid of pairs $(M, \Phi)$ 
of $R_W((u))$-modules $M$ which are fpqc locally on $R$ free
of rank $d$, together with a  $R_W((u))$-linear isomomorphism
$$
\Phi : \phi^*M:=R_W((u))\otimes_{\phi, R_W((u))}M\xrightarrow{ \ } M\ .
$$
Again it is easy to see that $\RR_d$ gives a stack for the fpqc topology.
Write   
$$
\r:\Cc_d\rightarrow \RR_d \ ;\quad   (\Mfr, \Phi)\mapsto (\Mfr[1/u], \Phi) 
$$
for the forgetful morphism. 

\medskip

 Fix an integer $m\geq 0$. Let us consider the stack $\Cc_{m,d}$ such that $\Cc_{m,d}(R)$ 
is the groupoid of $R_W[[u]]$-$\Phi$-modules $(\Mfr, \Phi)$ as above that
satisfy the additional hypothesis
\begin{equation}
u^m\Mfr\subset \Phi(\phi^*{\Mfr})\subset u^{-m}\Mfr \ .
\end{equation}
Once again, $\Cc_{m,d}$ gives a stack for the fpqc topology.
The natural morphism $\Cc_{m,d}\rightarrow \Cc_d$ is a representable 
closed immersion.

If $d$ is fixed we will often write $\Cc$, $\RR$, $\Cc_m$ instead of $\Cc_d$, $\RR_d$ and
$\Cc_{m,d}$.

\subsection{} For simplicity, we will set $G={\rm Res}_{W/\Z_p}GL_d$. 
Set 
\begin{equation*}
\begin{aligned}
LG(R)&:=GL_d(R_W((u)))\ ,\\
L^{+}G(R)&:=GL_d(R_W[[u]])\ ,\\
LG^{\leq m}(R)&:=\{A\in GL_d(R_W((u)))\ |\  A, A^{-1}\in u^{-m} M_d(R_W[[u]])\  \}\ .
\end{aligned}
\end{equation*}
Hence $L^{+}G=LG^{\leq 0}$. Note that the functor
$$
R\mapsto LG^{\leq m}(R)
$$
is represented by a scheme $LG^{\leq m}$ (which is infinite dimensional). Let $(\Mfr, \Phi)\in \Cc_m(R)$ such that $\Mfr$ is a free $R_W[[u]]$-module. 
By picking a $R_W[[u]]$-basis of $\Mfr$, we can write $\Phi$ as  multiplication by   $A\in LG^{\leq m}(R )$. Changing the basis 
by $g\in GL_d(R_W[[u]])$ amounts to changing $A$ to $g^{-1}\cdot A\cdot \phi(g)$. Therefore, we can write
\begin{equation}
\Cc_{m,d}=[  LG^{\leq m}/_\phi\, L^+G  ]
\end{equation}
where the quotient $/_\phi$ is via the right action of $L^+G(R )=GL_d(R_W[[u]])$ by $\phi$-conjugation
by $ A\star g=g^{-1}\cdot A\cdot \phi(g)$. 

Similarly, we can write
\begin{equation}
\Cc_{ d}=[  LG /_\phi\, L^+G  ]\ ,\quad \RR_{ d}=[  LG /_\phi\, LG  ].
\end{equation}
In fact, we can consider the fpqc stack $\ti \Cc_d$ defined as follows: $\ti\Cc_d(R)$ is
the groupoid of pairs $((\Mfr, \Phi), \alpha)$ of $R_W[[u]]$-$\Phi$-modules $(\Mfr, \Phi)$
together with an $R_W[[u]]$-module isomorphism 
$$
\alpha: R_W[[u]]^d\xrightarrow{\sim} \Mfr.
$$
The stack $\ti\Cc_d$ is represented by the ind-scheme $LG$
and the forgetful morphism
$$
\pi : \ti\Cc_d\rightarrow \Cc_d
$$
is a $L^+G$-torsor.

\subsubsection{} 
Denote by $\F_G=LG/L^+G$ the affine Grassmannian 
of $R_W[[u]]$-``lattices" in $R_W((u))^d$.
(Here by $R_W[[u]]$-lattice we mean a locally on $R$ free
$R_W[[u]]$-submodule $L$ of $R_W((u))^d$ such that $L\otimes_{R[[u]]}R((u))=R_W((u))^d$.)
The fpqc quotient $\F_G=LG/L^+G$ is represented by an ind-scheme
which is ind-projective over $\Z_p$. For $m\geq 0$, let $\F^{\leq m}_G$ be the projective subscheme 
of $\F_G$ parametrizing $R_W[[u]]$-lattices $L$ 
 $$
u^mR_W[[u]]^d\subset L\subset u^{-m}R_W[[u]]^d.
$$
(This is a finite union of Schubert varieties in the affine Grassmannian.)  Set 
$U_0(R)=L^+G (R)=GL_d(R_W[[u]])$ and define for $n\geq 1$ the principal congruence subgroup $U_n$ of level $n$ by   $U_n(R)=I+u^n\cdot M_d(R_W[[u]])$.  The 
subgroup scheme $U_n$ is normal in $L^+G$
and the quotient $L^+G/U_n$ is represented by the smooth finite type
group scheme $\Gg_n$ given by 
the Weil restriction of $GL_d$ from $W[[u]]/(u^n)$ to $\Z_p$
(so that $\Gg_n(R)=GL_d(R_W[[u]]/(u^n)$). Note that under the action of $L^+G$ on $\F^{\leq m}_G$ 
the subgroup $U_{2m}$ acts trivially, and hence the action factors through $\Gg_{2m}$. 

\begin{thm}\label{thm1}
a) For $m\geq 1$,  ${\cal C}_m=[LG^{\leq m}/_\phi L^+G]$ is an Artin stack of finite type over $\Z_p$.
We can write ${\cal C}$  as a direct $2$-limit
$$
{\cal C}=\lim_{\stackrel {\rightarrow}m } {\cal C}_m
$$
and so ${\cal C}$ is an ``ind-Artin stack of ind-finite type over $\Z_p$".

b) There is a formally smooth morphism 
$$
q: {\cal C}\rightarrow [L^+G \backslash \F_G]=[L^+G \backslash LG/L^+G]\ .
$$
In fact, $q$ is given as the limit of formally smooth morphisms 
$$
q_m: {\cal C}_m\rightarrow  [L^+G \backslash \F^{\leq m}_G]=[L^+G \backslash LG^{\leq m}/L^+G]\ .
$$ 
The composition of  $q_m$ with the natural morphism  $ [L^+G \backslash \F^{\leq m}_G]\rightarrow  [\Gg_{2m} \backslash \F^{\leq m}_G]$ is a smooth morphism of Artin stacks of finite type,
$$
\bar q_m: {\cal C}_m\rightarrow  [\Gg_{2m} \backslash \F^{\leq m}_G]\ .
$$
The relative dimension of $\bar q_m$ is equal to $2md^2$.  \end{thm}

\begin{proof}
The group $LG(R)$ is a topological group with topology described
 by
the neighborhoods $U_n$ of the identity $I=I_d$. We have
$$
LG(R)=\bigcup _{m\geq 0} LG^{\leq m}(R )\ .
$$
Suppose now that $A$ is in $LG^{\leq m}(R )$. 
For any integer $n\geq 0$ and for all $A'$ in the neighborhood (coset)
$$
\{A'\ |\    A'\cdot A^{-1} \in U_n(R)\}=    U_n(R)\, A
$$
of $A$, we have $A'\in LG^{\leq m}(R )$ also.

\begin{prop}\label{pr1} Suppose $n>2m/(p-1)$. 

1) For each $g\in U_n(R)$, $A\in LG^{\leq m}(R )$, we can write 
$ g^{-1}\cdot A\cdot \phi(g) = H^{-1}\cdot A $ with a unique 
$H=H(g,A)\in U_n(R)$. 
\begin{comment}The map $g\mapsto H(g, A)$ satisfies
\begin{equation}\label{coc}
H(gg', A)=H(g ,  H(g', A)\cdot A)\ .
\end{equation}
\end{comment}

2) Conversely, for each $A\in LG^{\leq m}(R ) $ and $h\in U_n(R)$,
there is a unique $g\in U_n(R)$ such that $A\star g=g^{-1}\cdot A\cdot \phi(g)= h^{-1}\cdot A $.
\footnote {The fact  that two elements $A$ and $A'$ in $GL_d(\bar\Ff_p((u)))$ 
which are $u$-adically close, are $\phi$-conjugate is also used  by Caruso in \cite{Ca}. The analogous fact for classical Dieudonn\'e modules is also true.  }
\end{prop}

\begin{proof} Let us first prove (1). Write   $g^{-1} =I+u^n X$ with $X\in M_d(R_W[[u]])$. Then $\phi(g) =I+u^{pn} Y$,
with $Y\in M_d(R_W[[u]])$. Now
$$
 g^{-1} \cdot A\cdot \phi(g) \cdot A^{-1}= (I+u^n X)\cdot A\cdot (I+u^{pn}Y)\cdot A^{-1}=
 $$
 $$=(I+u^nX)\cdot (I+u^{pn} AYA^{-1}).
$$
Observe that $A YA^{-1}\in u^{-2m}M_d(R_W[[u]])$ and $pn-2m>n$. Hence, for
$$
H^{-1}=(I+u^nX)\cdot (I+u^{pn} AYA^{-1})
$$ 
we obtain $g\cdot A\cdot \phi(g)^{-1}=H^{-1} \cdot A$.
The element $H$ is uniquely determined from $g$ and $A$
by   $g\cdot A\cdot \phi(g)^{-1}=H^{-1} \cdot A$. %The identity 
%(\ref{coc}) follows immediately.

The statement (2) is little trickier. First we show that if such a $g$ exists
it is uniquely determined by $h$ and $A$.
It is enough to assume $g\cdot A\cdot \phi(g)^{-1}=A$ 
with $g\in U_n(R)$ and $A\in LG^{\leq m}(R)$ and show $g=1$.
Write   $g=I+u^nX$, $\phi(g) =I+u^{pn}\phi(X)$. We have
$$
(I+u^nX)\cdot A=A\cdot (I+u^{pn}\phi(X))
$$
which gives
$$
u^n X\cdot A=u^{pn} A\cdot \phi(X)\ ,
$$
 i.e.,
$$
X_0+X_1u+X_2u^2+\cdots =u^{(p-1)n}  A\cdot (X_0+X_1u^p+X_2u^{2p}+\cdots )A^{-1}
$$
Note that $A\cdot X_i\cdot A^{-1}\in u^{-2m}M_d(R[[u]])$. Since $(p-1)n-2m>0$, we obtain $X_0=0$
which implies $g\in U_{n+1}(R)$. An induction finishes the proof of uniqueness. 

Now we will show that such a $g$ exists.
Let $A'=h^{-1}A$. Set $A_0=A$, $h_0=h$ and  
define $h_i$ and $A_i$ inductively by 
\begin{equation}\label{ind}
A_i=h_{i-1}^{-1}\cdot A_{i-1}\cdot \phi(h_{i-1})\ ,\quad  A'=h_i^{-1}\cdot A_i\ .
\end{equation}
Set $\kappa(i):=p^in-2(1+p+\cdots +p^{i-1})m$
with $\kappa(0)=n$. Note that under our assumption 
$n>2m/(p-1)$, the function
$\kappa(i)$  strictly increases with $i\geq 0$.

The existence of $g$ will now follow from 

\begin{lemma}
1) We have $h_i\in U_{\kappa(i)}(R)$ and so $\lim_{i\rightarrow\infty} h_i=I$.

2) Let $g_i=\prod_{j=0}^ih_j$. Then the limit $g=\lim_{i\rightarrow\infty}g_i $ exists and belongs to $U_n(R)$.
\end{lemma}
\begin{proof}
We will prove (1) by induction. It is true by our hypothesis when 
$i=0$. The equalities (\ref{ind}) imply
$$
h_i=A'\cdot \phi(h_{i-1})\cdot A'^{-1}.
$$
By the induction hypothesis $h_{i-1}\in U_{\kappa(i-1)}(R)$
so
$$
\phi(h_{i-1})=I+\phi(u^{\kappa(i-1)}X)=I+u^{p\cdot\kappa(i-1)}\phi(X) 
$$
with $X\in M_d(R_W[[u]])$. Since $A'\in LG^{\leq m}(R)$,
we have 
$$
A'\cdot \phi(X)\cdot A'^{-1}\in u^{-2m}M_d(R_W[[u]]),
$$
and so
$$
h_i=A'\cdot \phi(h_{i-1})\cdot A'^{-1}=I+u^{p\cdot\kappa(i-1)-2m}Y
$$
with $Y\in M_d(R_W[[u]])$. Since $\kappa(i)=p\cdot \kappa(i-1)-2m$ this completes the proof of (1). Part (2) now follows immediately  since
from part (1)
$$
g_i=\prod_{j=0}^ih_j=(I+u^{\kappa(0)}X_0)\cdot (I+u^{\kappa(1)}  X_1)\cdot\cdots\cdot (I+u^{\kappa(i)}X_i)
$$
with $X_j\in M_d(R_W[[u]])$ and $i\mapsto \kappa(i)$ is strictly increasing.
\end{proof}

\smallskip

Now $h_i^{-1}A_i=g_{i}^{-1}\cdot A\cdot \phi(g_{i-1})$, hence passing to 
the limit, we obtain
$
g^{-1}\cdot A\cdot \phi(g)=A'=h^{-1}\cdot A\
$
as desired.
\end{proof}

 \begin{Remarknumb} \label{drinfeld}{\rm Let $M$ be a  $R((u))$-module 
and let $R\rightarrow R'$ be a flat extension such that 
$M'=M\hat\otimes_RR'\simeq R'((u))^d$. The module $M$ has 
a natural topology as a Tate $R$-module (see  [Dr]). 
The $R$-lattices  of $M$ (i.e  $R$-modules $L$
which are open and  such that $L/U$ is finitely generated for every open submodule $U\subset L$) give a basis
of open neighborhoods of $0$.
Multiplication by $u$ on $M$ is topologically nilpotent;
i.e given any two $R$-lattices $L$, $L'$, there is $N\geq 0$ such that
$u^N\cdot L\subset L'$.

Then $\Gg_M:={\rm Aut}_{R((u))}(M)$ has a natural structure of a topological group. 
To obtain a basis of neighborhoods of the identity $I$ we take a lattice $L$ and for $n>>0$ we consider
$$
U_n(L)=\{g\in \Gg_M\ |\  g(L)\subset L,\ \ g_{|L}\equiv I\ {\rm mod}\ u^nL\}.
$$
We can then show:
\smallskip

{\sl Given two isomorphisms $A, A': \phi^*M\xrightarrow{\sim} M$,
there exists an open neighborhood $U$ of the identity  in the topological group $\Gg_M$ such that 
if $A'\cdot A^{-1}\in U$ then $A$, $A'$ are $\phi$-conjugate
by a uniquely determined element of $\Gg_M$.}
\smallskip

The argument is similar as above: Suppose that $A: \phi^*M\xrightarrow{\sim} M$ is an $R((u))$-isomorphism and for $h\in \Gg_M$ define $h_0=h$ and 
inductively
$$
h_i=A\cdot \phi^*(h_{i-1})\cdot  A^{-1}.
$$
The result follows from the statement:
There is an open neighborhood $U$ of $I$ such that 
for $h\in U$, we have $\lim_{i\rightarrow \infty} h_i=I$
and the limit $\tilde h=\lim_{i\rightarrow \infty} \prod_{0\leq j\leq i}h_j$
exists. Indeed, the arguments above show that this is true when $M$ is a free $R((u))$-module. In general, let $R\rightarrow R'$ be a flat 
homomorphism such that $M'=M\hat\otimes_RR'\simeq R'((u))^d$.
Consider $h'_i=h_i\otimes 1\in \Gg_{M'}$ and let $L$ be a lattice in $M$. Then $L'=L\hat\otimes_RR'$ is a lattice in $M'$. By the above, there is $n$ such that when $h\in U_n(L)$ (and hence $h'=h\otimes 1\in U_n(L')$) we have
$$
(h_i(x)-x)\otimes 1=h'_i(x')-x'\in u^{\kappa(i)}L'
$$
for a strictly increasing sequence $\kappa(i)$.
Since $u^{\kappa(i)}L'\cap M=u^{\kappa(i)}L$ this shows 
the result.}
\end{Remarknumb}

We now continue with the proof of Theorem \ref{thm1} (a).
Recall
$$
\Cc_m=[LG^{\leq m}/_\phi\, L^+G]
$$
Let $n>2m/(p-1)$. Recall the normal subgroup $U_n$ of $L^+G$ and its smooth finite type quotient $\Gg_n$. Consider the quotient
stack $[LG^{\leq m}/_\phi\, U_n]$. Proposition \ref{pr1} implies that
$[LG^{\leq m}/_\phi\, U_n]$ coincides with  the 
quotient $X_{n, d}^{\leq m}:=[LG^{\leq m}/  U_n]$ by the free translation 
action of $U_n$ on $LG^{\leq m}$. The quotient $X_{n,d}^{\leq m}$ 
is represented by a scheme of finite type over $\Z_p$.
This can be seen as follows. 

Recall that the quotient $X_{0, d}^{(m)}=[LG^{\leq m}/  L^+G]$
is represented by the closed subscheme $\F_G^{\leq m}$ of the affine Grassmannian $\F_G=LG/L^+G$ that parametrizes lattices $L$ such that
$
u^mL_0\subset L\subset u^{-m}L_0 
$,
where $L_0=R_W[[u]]^d$.
The natural map
$$
p^{\leq m}: X_{n,d}^{\leq m}\rightarrow X_{0, d}^{\leq m}=\F^{\leq m}_G 
$$
is represented by the $\Gg_n$-torsor
that parametrizes pairs $(L, \alpha)$ where $L$ 
is a lattice as above and $\alpha: L_0/u^n L_0\xrightarrow{\sim}
L/u^nL$ is an $R_W[[u]]/(u^n)$-isomorphism.

Combining the above, we  now see that 
\begin{equation}\label{quo2}
\Cc_m\simeq [X^{\leq m}_{n,d}/_\phi\, \Gg_n]\ ,
\end{equation}
where the quotient is for the action of the smooth group scheme
$\Gg_n$ on $X^{\leq m}_{n, d}$ which is induced by $\phi$-conjugation.
This (right) action of $\Gg_n$ on $X^{\leq m}_{n, d}$  can be explicitly described as follows:
Let $\gamma\in \Gg_n(R)$ which  we can lift to
$g\in L^+G(R)=GL_d(R_W[[u]])$ and consider the point $x=(L, \alpha)\in  X^{\leq m}_{n, d}(R)$
given through the matrix $A$ by $L= L_0\cdot A$, $\alpha=A\ \hbox{\rm mod}\ u^n$.
 (Here the elements of $L_0$ are viewed as row vectors.)
 Then $ x\star \gamma$ is the point of $X^{\leq m}_{n, d}$ that corresponds to the matrix 
 $g^{-1}\cdot A\cdot \phi(g) \in LG^{\leq m}(R)$.
 Observe that if $n'>n>2m/(p-1)$, the natural morphism $X^{\leq m}_{n', d}\rightarrow X^{\leq m}_{n,d}$ induces an isomorphism
\begin{equation}\label{quo3}
 [X^{\leq m}_{n',d}/_\phi\, \Gg_{n'}]\xrightarrow{ \sim }[X^{\leq m}_{n,d}/_\phi\, \Gg_n\ ]\, .
\end{equation} 
It follows from (\ref{quo2}) and the above that $\Cc_m$ is an algebraic
 (Artin) stack of finite type over $\Z_p$ of   dimension 
 equal to the   dimension of the scheme $\F^{\leq m}_G$.
 It is clear that we can represent $\Cc$ as the $2$-limit of 
 the algebraic stacks $\Cc_m$ 
  and so the rest of (a) follows.
 
For part (b) observe that the quotient description of $\Cc_m$ 
implies the existence of 
$$
q_m: \Cc_m=[LG^{\leq m}/_\phi\, L^+G]\rightarrow [L^+G\backslash \F^{\leq m}_G]=[ \F^{\leq m}_G/L^+G]
$$
(here in the last quotient $g\in L^+G$ acts by $L\cdot g=g^{-1}L$).
This descends the quotient morphism
$$
LG^{\leq m}\rightarrow LG^{\leq m}/L^+G=\F^{\leq m}_G\ .
$$
Now for $n\geq 2m$, $U_n$ acts trivially on $\F^{\leq m}_G$
and the action of $L^+G$ on $\F^{\leq m}_G$ factors through the quotient $\Gg_n$. Hence composing $q_m$ with the quotient morphism by $U_n$, we obtain a morphism  
$q'_m: \Cc_m\rightarrow  [\F^{\leq m}_G/\Gg_n]$. When $n>2m/(p-1)$, the morphism $q'_m$ is given by taking the quotient 
$$
[X^{\leq m}_{n,d}/_\phi\, \Gg_n]\rightarrow [\F^{\leq m}_G/\Gg_n] 
$$
of the smooth torsor $X^{\leq m}_{n,d}\rightarrow \F^{\leq m}_G$
and hence is smooth. It follows that $q_m$ itself is formally smooth.
Also the morphism $\bar q_m$ of the statement of part (b), which is 
given as a composition of
$q'_m$ with $[ \F^{\leq m}_G/\Gg_n] \rightarrow  [\F^{\leq m}_G/\Gg_{2m}] $,
is also smooth. A straightforward dimension count
now gives that the relative dimension of $\bar q_m$ is equal to
the (relative) dimension of $\Gg_{2m}$ over $\Z_p$; this is equal to $2md^2$.
\end{proof}
 
 \begin{comment}
 \subsection{}
 As an example, suppose that $m$ is small compared to $p-1$ so that we can take $n=1$.
Then we have $GL_d(R[[u]])/U_n(R)=GL_d(R)$ and
$$
\Cc_m=[GL_d\backslash _{\phi\ } X^{(m)}_{1,d}]\ .
$$
In this case, $GL_d(R[[u]])\to GL_d(R[[u]])/U_n(R)=GL_d(R)$ splits using the natural 
$GL_d(R)\subset GL_d(R[[u]])$. Hence, 
 the $\phi$-conjugation  amounts to standard conjugation 
 by $g\in GL_d(R)\subset GL_d(R[[u]])$.
We then have
$$
\Cc_m=[GL_d\backslash_{\rm conj}\   X^{(m)}_{1,d}]\rightarrow X^{(m)}_{0,d}\subset \F_d\ .
$$
\end{comment}

\subsection{}
We consider now some properties of $\RR$,  $\theta_m:\Cc_m\rightarrow \RR$ and $\r: \Cc\rightarrow \RR$. Recall that, for each $R_W((u))$-module $M$ which is fpqc locally on $S=\Spec(R)$
free of rank $d$, we have the (twisted) affine Grassmannian $Gr_M\rightarrow S$ whose $A$-points for an $R$-algebra $A$ 
are given by 
$A_W[[u]]$-lattices $\Mfr$ of $M_A=M\hat\otimes_RA$. By [Dr] (Theorem 3.8 and  Remark (b) below it),
$Gr_M$ is represented by an ind-algebraic space which is ind-proper and of ind-finite presentation over $S$.

\begin{thm}\label{grass}
a) For each $S=\Spec(R)\rightarrow \RR$
which corresponds to a $R_W((u))$-$\Phi$-module $(M, \Phi)$,
the fiber products
$$
\r \times_\RR S: \Cc \times_\RR S\rightarrow S \ ,\quad \r_m \times_\RR S: \Cc_m \times_\RR S\rightarrow S,
$$
are represented  by  the 
(twisted) affine Grassmannian $Gr_M\rightarrow S$,
resp. by a proper algebraic subspace of 
$Gr_M\rightarrow S$.

b)  The diagonal morphism 
$\delta: \RR\rightarrow \RR\times_{\Z_p} \RR$ is representable and of finite presentation.
\end{thm}

\begin{cor}\label{grasscor}
a)   
$ 
\r : \Cc  \rightarrow \RR
$ 
is ind-representable and ind-proper. 

b) 
$ 
\r_m: \Cc_m\rightarrow \RR
$ 
is representable, proper and of finite presentation. 
\end{cor}
  
\begin{proof}
Part (a). The first part of the statement regarding 
$\r\times_{\RR} S:\Cc \times_\RR S\rightarrow S$ 
follows from the definition. 
Note here that we do not necessarily know that $M$ contains a free $R_W[[u]]$-lattice.
However, there is a flat 
base change $R\rightarrow R'$ 
such that  $M'=M\hat\otimes_R R'$ is $R'_W((u))$-free.
Then there is a (free) $R'_W[[u]]$-lattice $\Mfr'_0$
in $M'$. We will now show the second part of the statement.
Let $\delta$ be the smallest integer for which
\begin{equation}\label{inclu2}
u^\delta \Mfr'_0\subset \Phi(\phi^*\Mfr'_0)\subset u^{-\delta}\Mfr'_0.
 \end{equation}
Set $S'=\Spec(R')$ with $R'$ as above.
Suppose $T=\Spec(A)$ is an $S$-scheme and set $A'=A\otimes_RR'$, $T'=T\times_SS'=\Spec(A')$.
Let $\Mfr$ be an $A'_W[[u]]$-$\Phi$-lattice in $M\hat\otimes_RA'$
that corresponds to an object  of 
$(\Cc_m\times_\RR S')(T')$. For simplicity, set $\Mfr'_{0,A}=\Mfr'_0\hat\otimes_{R'} A'$.
Then
\begin{equation}\label{inclu}
u^m\Mfr\subset \Phi(\phi^*\Mfr)\subset u^{-m}\Mfr\ , \quad u^N\Mfr'_{0,A}\subset \Mfr\subset u^{-N}\Mfr'_{0,A} ,
\end{equation}
for some $N\geq 0$ (we can suppose that $N$ is the smallest integer with this property).
Applying $\Phi$ to the second chain of inclusions (\ref{inclu}) gives
$$
u^{pN}\Phi(\phi^*\Mfr'_{0,A})\subset \Phi(\phi^*\Mfr)\subset u^{-pN}\Phi(\phi^* \Mfr'_{0,A})
$$
and $pN$ is the smallest integer with this property. On the other hand, we have
$$
\Phi(\phi^*\Mfr)\subset u^{-m}\Mfr\subset u^{-m-N}\Mfr'_{0,A}\subset u^{-m-N-\delta}\Phi(\phi^*\Mfr'_{0,A}), \hbox{\rm \ and}
$$
$$
u^{m+N+\delta}\Phi(\phi^*\Mfr'_{0,A})\subset u^{m+N}\Mfr'_{0,A}\subset u^m\Mfr\subset \Phi(\phi^*\Mfr). \ \ \ \ \ \ \ \ \
$$
Combining these gives
$$
u^{N+m+\delta}\Phi(\phi^*\Mfr'_{0,A}) \subset \Phi(\phi^*\Mfr)\subset u^{-N-m-\delta}\Phi(\phi^*\Mfr'_{0,A}).
$$
This implies that
$pN\leq N+m+\delta$, i.e $N\leq (m+\delta)/(p-1)$, and so
$$
\quad u^{[\frac{m+\delta}{p-1}]}\Mfr'_{0,A} \subset \Mfr\subset u^{-[\frac{m+\delta}{p-1}]}\Mfr'_{0,A} .
$$
(This is essentially the same argument as in [Ki1] Prop. 2.1.7.)
By the above, and the definition [Dr] of the ind-structure on $Gr_M$, this implies that $\Cc_m\times_\RR S'$ is represented by a proper $S'$-scheme; therefore, by descent,  $\Cc_m\times_\RR S$ is
an $S$-proper   algebraic space.

Part (b). Suppose that $M$, $N$ are two $R_W((u))$-$\Phi$--modules
of rank $d$. Consider the functor on $R$-algebras
$$
A\mapsto {\rm Isom}_{\RR}(M, N)(A):={\rm Isom}_{\Phi, A_W((u))}(M_A, N_A),
$$
where for simplicity we write $M_A=M\hat\otimes_RA$,
$N_A=N\hat\otimes_RA$. We will show that
this is representable by a scheme of finite
presentation over $R$. This  implies then the statement in (b).
Using the existence of $R_W$-lattices in both $M$ and $N$ ([Dr], see Remark \ref{drinfeld}),
we can see that the functor that sends $R$
to the $R_W((u))$-linear isomorphisms $M\rightarrow N$
is represented by an ind-scheme.
It is not hard to see that ${\rm Isom}_{\RR}(M, N)$ is represented by a ind-closed ind-subscheme of this ind-scheme. 
To show that this is actually a scheme of finite
presentation we can employ an fpqc base change $R\rightarrow R'$
and assume that $M'=M_{R'}$, $N'=N_{R'}$ are given by  
 $A$, $B\in LG(R')$.  
By the definitions, the additional condition on the
$R'_W((u))$-linear isomorphism  $M'\rightarrow N'$
given by 
$g\in GL_d(R'_W((u)))$ that guarantees that it respects $\Phi$ 
is  
\begin{equation*}
A=g^{-1}\cdot B\cdot \phi(g), \quad \hbox{\rm or equivalently,}
\end{equation*}
\begin{equation}\label{id}
g =B\cdot \phi(g)\cdot A^{-1}\ .
\end{equation}
Suppose that $A$ and $B$ are in $LG^{\leq m}(R')$,
$LG^{\leq n}(R')$ respectively. Assume that $g$ is in $LG^{\leq s}(R')$ and $s=s(g)$ is the smallest integer
with that property. Then $\phi(g)$ belongs to $LG^{\leq ps}(R')$
and we can see that $ps$ is the smallest 
integer with this property.
The identities above now imply that
$\phi(g)$ is in $LG^{\leq s+m+n}(R')$.
Therefore $ps\leq s+m+n$ which gives
\begin{equation}\label{ineq}
 s\leq \frac{m+n}{p-1}\ .
\end{equation}
Let us write out (\ref{id}) explicitly
\begin{equation}\label{id2}
\sum_{i\geq -s} g_i u^i=B\cdot (\sum_{i\geq -s} g_i u^{pi})\cdot A^{-1}=\sum_{i\geq -s} u^{pi}\cdot (B\cdot g_i\cdot A^{-1})
\end{equation}
with $g_i\in M_d(R')$.  Now  
consider the matrix identity obtained by 
comparing the $u^a$ 
terms of both sides of (\ref{id2}) for $a>(m+n)/(p-1)$.
We see that this has the form
\begin{equation}
g_a=\sum_{i, k, l} B_l\cdot g_{i}\cdot A'_{k}
\end{equation}
with $pi+k+l=a$ and $i\geq -s$, $k\geq -m$, $l\geq -n$
and $A^{-1}=\sum_{k=-m}^\infty A'_k u^k$. Since 
$a>(m+n)/(p-1)$, these inequalities imply that $i<a$.
Therefore, all these matrix identities for $a>(m+n)/(p-1)$
amount to determining $g_a$ from $g_i$ for $i<a$.
 The result now follows.
\end{proof}
 
 \medskip

\begin{cor}\label{dia3}
There is a diagram
\begin{equation*}
\xymatrix{
{} & \ar[dl]_\r \Cc \ar[dr]^q & {}\\
\RR & {} &   \ \ [L^+G\backslash \F_G]\ ,
 }
\end{equation*}
where  the morphism
$q$ is formally smooth and the morphism
$\r$ is ind-representable and ind-proper.
\end{cor}

\subsubsection{}\label{schubert}
Recall that we set $G={\rm Res}_{W/\Z_p}GL_d$. Let $K_0$ be the fraction field of $W$ and
set $f=[K_0:\Q_p]=[k:\Ff_p]$. Then, after ordering the elements of ${\rm Gal}(K_0/\Q_p)$,
we can write
$$
G(K_0)=\prod_{i=1}^f GL_d(K_0), \quad (\F_G)_W=\prod_{i=1}^f (\F_{GL_d})_W \ .
$$
Let us set
$\nu=(\nu{(1)},\ldots , \nu{(f)})$ where
for each $i=1,\ldots , f$, $\nu{(i)}=(n_1(i),\ldots , n_d(i))$
is a collection of integers with $n_1(i)\geq n_2(i)\geq \cdots \geq n_d(i)$.
Let $F\subset K_0$ be the fixed field of the subgroup of ${\rm Gal}(K_0/\Q_p)$
that fixes $\nu$.
 
%Let 
%$$
%\nu: K_0^*\rightarrow G(K_0)=\prod_{i=1}^f GL_d(K_0)
%$$
%be a cocharacter that we can write $\nu=(\nu{(1)},\ldots , \nu{(f)})$.
Denote by $u^{\nu(i)}$ the diagonal matrix ${\rm diag}(u^{n_1(i)},\ldots , u^{n_d(i)})$
in $GL_d(W((u)))$
and set
$$
u^{\nu}= (u^{\nu(1)},\ldots , u^{\nu(f)}) \in LG(W).
$$
Suppose that $N={\rm max}\{|n_j(i)|\}$. Let $S^0_\nu$, resp. $S_\nu$, be the corresponding open, resp. projective, affine Schubert variety in $(\F^{\leq N}_G)_W\subset 
(\F_G)_W$ which is given as the image of $L^+G\cdot u^\nu\cdot L^+G$,
resp. the Zariski closure of that image. By descent, we see that this is defined over the integers $W'$ of the reflex field $F$. We set ${\cal C}_\nu$ for 
the Artin stack over $W'$ which is the inverse image of $S_\nu$ under $q$; this is a closed substack of $({\cal C}_N)_{W'}$. If $n_d(i)\geq 0$ for all $i$, then
for $S=\Spec(R)$, the groupoid ${\cal C}_\nu(S)$ is given by $R_W[[u]]$-$\Phi$-modules $(\Mfr,\Phi)$ of rank $d$ such that $\Phi(\phi^*\Mfr)\subset \Mfr$ and such that the action of $u$ on $Coker(\Phi)=\Mfr/\Phi(\phi^*\Mfr)$ has elementary divisors $u^{\nu'}$ with $\nu'$ which satisfies $\nu'(i)\leq \nu(i)$ in the usual ordering, for all $i$; ${\cal C}_\nu$ is a closed substack of ${\cal C}_{N}$. 

The obvious
version of Theorem \ref{thm1} holds for the stack ${\cal C}_\nu$; it is an Artin stack of finite type over $W'$ smoothly equivalent to the Schubert variety $S_\nu$ in the
affine Grassmannian $\F_G$.
Furthermore, restricting  $\r$
to ${\cal C}_\nu$ gives 
$$
\r_\nu: {\cal C}_\nu\rightarrow {\cal R}_{W'}\ , 
$$
which is representable, proper and of finite presentation. 
(Similarly, we can consider 
the inverse image ${\cal C}_\nu^0$ of $S^0_\nu$ under $q$ and the restriction 
 $\r_\nu: {\cal C}_\nu^0\rightarrow {\cal R}_{W'} $
which is therefore representable and of finite presentation.)
Summarizing, we obtain a diagram
\begin{equation}\label{dia2}
\xymatrix{  & \ar[dl]_{\r_\nu} \Cc_\nu \ar[dr]^q &  \\
\RR_{W'} &  &   \ \ [(L^+G)_{W'}\backslash S_\nu]\ ,}
 \end{equation}
where  the morphism
$q$ is  formally smooth and the morphism
$\r_\nu$ is  representable and  proper.

\subsection{} \label{genG}
Let us sketch how to generalize the above theory to reductive groups. Let $H$ be any reductive algebraic group scheme $H$ over $W$. Let us set $G={\rm Res}_{W/\Z_p}(H)$. Instead of $R_W[[u]]$-$\Phi$-modules
of rank $d$ we consider $H$-torsors $\cal T$ over $R_W[[u]]$ together with a $H$-isomorphism $\Phi: \phi^*_S ({\cal T}[1/u])\xrightarrow{\sim} {\cal T}[1/u]$
(here $S=\Spec(R)$). 
The corresponding fpqc stack $\Cc_G$ can be viewed as the quotient 
$[LG/_\phi L^+G]$. Similarly, we can consider the fpqc stack $\RR_G$ of
$H$-torsors $T$ over $R_W((u))$ (which are trivial fpqc locally
on $R$),  together with a $H$-isomorphism $\Phi: \phi^*_S (T)\xrightarrow{\sim} T$. 
The stack $\RR_G$ can be viewed as the quotient 
$[LG/_\phi\, LG]$. The obvious generalizations of Theorems \ref{thm1} and Corollary \ref{grasscor}  
hold in this situation. The proofs are extensions of 
the above proofs for $GL_d$ after using a faithful representation $H\hookrightarrow GL_N$.
For example, the ind-structure is modeled on the ind-scheme structure of $LG=\lim_{m} (LG\cap L{\rm Res}_{W/\Z_P}GL_N^{\leq m})$.
In particular, we again obtain a diagram
\begin{equation}\label{dia3}
\xymatrix{
  & \ar[dl]_\r \Cc_G \ar[dr]^q & \\
\RR &   &   \ \ [L^+G\backslash \F_G]\ ,}
\end{equation}
where  the morphism
$q$ is formally smooth and the morphism
$\r$ is ind-representable and ind-proper.

 Suppose we are given a
dominant coweight $\nu$ of $G$ defined over 
an unramified extension $F$ of $\Q_p$ with integers $W'$.
Denote by $u^\nu\in LG(F)=G(F((u)))$ the element given as the image of $u\in \Gm(F((u)))$
by the corresponding homomorphism
$
\Gm_F\rightarrow G_F
$.
As in \S \ref{schubert}, we can define 
$S^0_\nu$, resp. $S_\nu$, to be the corresponding open, resp. projective, affine Schubert variety in 
$(\F_G)_{W'}$ which is given as the image, resp. the closure of the image of $(L^+G)_{W'}\cdot u^\nu\cdot (L^+G)_{W'}$.

We set ${\cal C}_{G,\nu}$ for 
the Artin stack over $W'$ which can be defined by descent
as the inverse image of $S_\nu$ under $q$.
The obvious
version of Theorem \ref{thm1} holds for the stacks  ${\cal C}_{G, \nu}$; they are Artin stacks of finite type over $W'$ smoothly equivalent to Schubert varieties in the
affine Grassmannian $\F_G$.
Once again, we have a diagram 
\begin{equation}\label{dia5}
(\RR_G)_{W'} \xleftarrow{ \r_\nu } {\cal C}_{G,\nu}  \xrightarrow{q_\nu\ } [(L^+G)_{W'}\backslash S_\nu]\ , 
\end{equation}
where  the morphism
$q_\nu$ is  formally smooth and the morphism
$\r_\nu$ is  representable and  proper.

    \bigskip
    \bigskip

\section{ $p$-adic models and local models}\label{locmod}

In this section, we define the stacks $\Cc_{d,h,K}$ 
and prove Theorem \ref{thmintro} of the introduction.
Let $K$ a finite extension of $\Q_p$ with residue field $k$ and ramification index $e$.
Choose a uniformizer $\pi$ of $K$ with Eisenstein polynomial $E(u)$ 
over $K_0={\rm Fr}(W)$. Then we can write $\O_K\simeq W[[u]]/(E(u))$.
Fix $h\geq 1$ (the ``height"). 
For $a\geq 1$, note that $u^{ea}$ vanishes in $W_a[[u]]/(E(u))$
where $W_a=W/p^aW$. Hence, if $R$ is a 
$\Z/p^a\Z$-algebra, we have
\begin{equation}\label{inclusion}
 u^{eah}R_W[[u]]\subset E(u)^hR_W[[u]].
\end{equation}

Now consider the category ${\rm Nil}_p$  of schemes $S$ such that 
$p^b\cdot \O_S=0$ for some $b\geq 1$. Such schemes can be viewed
as formal schemes over $\Z_p$.
We will call set-valued functors on ${\rm Nil}_p$ 
which satisfy descent for the fpqc topology ``formal spaces".
A formal scheme $X$  over ${\rm Spf}(\Z_p)$  
gives a formal space by sending $S$ in ${\rm Nil}_p$
to the set of formal scheme morphisms $S\rightarrow X$
over ${\rm Spf}(\Z_p)$. Also if $\cal S$ is a fpqc stack over
$\Z_p$, we can consider the restriction $\hat{\cal S}$ to a
groupoid over the category ${\rm Nil}_p$; we can think of the ``formal stack" $\hat{\cal S}$ as
``the formal completion of $\cal S$ along its fiber over $p$".

\subsection{} \label{local} Consider the functor $M_{d, h, K}$ on schemes
over $\Z_p$
that associates to $S=\Spec(R)$ the set of $R_W[u]$-submodules
\begin{equation}
 {\cal E}\subset (R_W[u]/(E(u)^h))^d \ ,
 \end{equation}
such that both $\cal E$ and the quotient
$(R_W[u]/(E(u)^h))^d/{\cal E}$ are $R_W$-projective
with rank locally constant on $\Spec(R)$.
This functor is represented by a projective scheme 
  over $\Z_p$ (a disjoint sum of closed subschemes of Grassmannians), which we will also denote by $M_{h, K}$.
 Once again, here and in what follows we will omit the subscript $d$ 
 from the notation. In fact, if in addition $h=1$, we will 
 also omit $h$ from the notation and simply write $M_K$.
 The group scheme 
 $$
 {\rm Res}_{(W[u]/(E(u)^h))/\Z_p}GL_d
 $$
 over $\Z_p$ acts on $M_{h,K}$.
 
 Suppose that $p^a\cdot R=0$ for $a\geq 1$. Then $M_{h, K}(R)$
 is in bijection with the set of $R_W[[u]]$-modules $\L$ with
 $$
 u^{eah}\cdot R_W[[u]]^d\subset E(u)^h R_W[[u]]^d\subset \L\subset R_W[[u]]^d
$$
which are, locally on $R$, free over $R_W[[u]]$. This gives a functorial injection, 
\begin{equation}\label{grassincl}
M_{h,K}(R)\hookrightarrow \F_G(R)
\end{equation}
and it implies that we can view the formal completion $\hat M_{h,K}$
of $M_{h, K}$ along its fiber over $p$ as a subspace of the formal space $\hat\F_G$
defined by the affine Grassmannian $\F_G$.

\subsection{}
We now define a groupoid $\Cc_{d,h,K}$ over $\Z_p$-schemes
as follows.  Let $R$ be a $\Z_p$-algebra. 
Then
  $\Cc_{d,h,K}(R)$ is given by
 pairs $(\Mfr, \Phi)$
of an $R_W[[u]]$-module $\Mfr$ which is, locally fpqc on $\Spec(R)$, free 
of rank $d$ and a 
$R_W((u))$-module isomorphism 
\begin{equation}
\Phi : \phi^*{\frak M}[1/u]\xrightarrow{\sim} {\Mfr}[1/u]  ,
 \end{equation} 
 such that  $E(u)^h\Mfr\subset \Phi(\phi^*\Mfr)\subset \Mfr$.
We can see that the groupoid $\Cc_{d,h,K}$ is an fpqc stack. 
In what follows, we will omit $d$ from the notation and write
$\Cc_{h, K}$.

We can consider the formal $p$-adic completion $\widehat\Cc_{h, K}$ of $\Cc_{h, K}$.
This is a fpqc stack 
over ${\rm Nil}_p$ defined by considering $\Cc_{h, K}(R)$ as above
for  $\Z_p$-algebras $R$ in which $p$ is nilpotent. 
We can write   $\hat\Cc_{h, K}$ as a $2$-limit
$$
\widehat\Cc_{h, K}:=\lim_{\stackrel \rightarrow a} \Cc^a_{h,K}\ 
$$
where $\Cc^{a}_{ h,K}:=\widehat\Cc_{h,K}\times_{\Z_p}{\Z/p^a\Z}$
is the reduction modulo $p^a$.
Using (\ref{inclusion})
we can see that $\Cc^{a}_{ h,K}$ is a (closed) substack of  $\Cc_{eah}\times_{\Z_p}\Z/p^a\Z$.
Now suppose that $\mathcal X$ is a formal scheme over ${\rm Spf}(\Z_p)$ such that 
$p\mathcal O_{\mathcal X}$ is an ideal of definition.
Then, for each $a\geq 1$, ${\mathcal X}\times_{\Z_p}\Z/p^a\Z$ is a scheme over $\Z/p^a\Z$.
We set
$$
\widehat\Cc_{h, K}({\mathcal X}):=\lim_{\stackrel \leftarrow a} \widehat\Cc_{h,K}({\mathcal X}\times_{\Z_p}\Z/p^a\Z)\ .
$$  
This allows us to extend $\widehat\Cc_{h, K}$ to a groupoid over the category of adic formal schemes  over ${\rm Spf}(\Z_p)$.
Suppose that $R$ is a Noetherian $p$-adic ring (i.e a Noetherian $\Z_p$-algebra which is $p$-adically complete and separated,  $R=\ds{\lim_{\leftarrow {a}}  R/p^aR}$),
and set ${\mathcal X}={\rm Spf}(R)$.
Set
$$
R_W\{\{u\}\}=\lim_{\stackrel \leftarrow {a}} \left[(R_W/p^aR_W)((u))\right]=\left\{\sum_{i=-\infty}^{+\infty}a_iu^i\ |\ a_i\in R_W, \lim_{i\rightarrow -\infty}a_i=0\right\}.
$$
We can see that the objects of 
  $\widehat\Cc_{h,K}({\rm Spf}(R))$ are given by
 pairs $(\Mfr, \Phi)$
of an $R_W[[u]]$-module $\Mfr$ which is locally 
$R_W[[u]]$-free of rank $d$ and a 
$R_W\{\{u\}\}$-module isomorphism 
\begin{equation}
\Phi : \phi^*{\frak M}\otimes_{R_W[[u]]}R_W\{\{u\}\} \xrightarrow{\sim} {\Mfr}\otimes_{R_W[[u]]}R_W\{\{u\}\}  ,
 \end{equation} 
 such that  
 \begin{equation}\label{inclusionMfr}
 E(u)^h\Mfr\subset \Phi(\phi^*\Mfr)\subset \Mfr\ .
 \end{equation}
 (Note that $E(u)$ is a unit in $W\{\{u\}\}$ and so therefore also in 
  $R_W\{\{u\}\}$.)

\subsubsection{} For $\Spec(R)$ in ${\rm Nil}_p$ now set
$$
LG^{h,K}(R)=\{A\in M_d(R_W[[u]])\  |\ A^{-1}\in E(u)^{-h}\cdot M_d(R_W[[u]])\subset M_d(R_W((u))) \}.
$$
This defines a  functor on ${\rm Nil}_p$.
As before, we can write
$$
\widehat\Cc_{h,K}=[LG^{h,K}/_\phi\, L^+G].
$$
where (by abusing notation) we also 
denote by $L^+G$  the formal $p$-adic completion 
of $L^+G$.
The map
$$
A\mapsto A\cdot R_W[[u]]^d\subset R_W((u))^d
$$ 
gives a morphism of formal stacks,
$$
q_{ h, K} : \widehat\Cc_{ h,K}=[LG^{h,K}/_\phi\, L^+G] \rightarrow [L^+G\backslash \widehat M_{ h,K}].
$$
 
 \subsubsection{} \label{3b2} Using Proposition
\ref{pr1} and the above, we see that if $n(a)> eah/(p-1)$ then 
$$
[LG^{h, K}/_\phi U_{n(a)}]_{\Z/p^a\Z}\simeq [LG^{h, K}/ U_{n(a)}]_{\Z/p^a\Z}.
$$
As in the proof of Theorem \ref{thm1} we can see that the quotient stack $[LG^{h, K}/ U_{n(a)}]_{\Z/p^a\Z}$ is represented by 
a torsor $(X^{h,K}_{n(a),d})_{\Z/p^a\Z}$ for the group scheme $(L^+G/U_{n(a)})_{\Z/p^a\Z}=(\Gg_{n(a)})_{\Z/p^a\Z}$
over   $(M_{h, K})_{\Z/p^a\Z}$. Similarly to that proof, we can conclude
$$
\Cc^a_{h,K}\simeq [(X^{h,K}_{n(a),d})_{\Z/p^a\Z}/_\phi\, (\Gg_{n(a)})_{\Z/p^a\Z}]\ ,
$$
and that the morphism
$$
q^a_{h,K}:   \Cc^a_{h,K}\rightarrow [L^+G\backslash   M_{h,K}]_{\Z/p^a\Z}
$$
is formally smooth. Hence, the morphism 
between the formal stacks
\begin{equation}\label{qhK}
q_{h,K}:   \widehat\Cc_{h,K}\rightarrow [L^+G\backslash   \widehat M_{h,K} ]
\end{equation}
is also formally smooth. 

Recall the definition of the stack $\mathcal R=\mathcal R_d$ over $\Z_p$-schemes, cf. 
\S \ref{2defofstacks}, and the notations  $\mathcal R^a=\RR\times_{\Z_p}\Z/p^a\Z $,
$\Cc_{h, K}^a=\Cc_{h, K}\times_{\Z_p}\Z/p^a\Z$. By sending an object $(\mathfrak M, \Phi)$ to $(\mathfrak M[1/u], \Phi)$, 
we obtain a morphism of stacks $\r: \Cc_{h, K}\rightarrow \mathcal R$. Note that the morphism 
$\widehat\r: \widehat\Cc_{h, K}\rightarrow \widehat\RR$ on formal
 completions is obtained by passing to the limit on the morphisms 
$\r^a:\Cc^a_{h, K}\rightarrow \RR^a$ which arise by restricting the morphisms 
$$
\r_{eah}\times_{\Z_p}\Z/p^a\Z : \Cc_{eah}\times_{\Z_p}\Z/p^a\Z \xrightarrow { \ } \RR^a
$$ 
to the closed substacks $\Cc^a_{h, K}$. This together with Corollary \ref{grasscor}
implies that the morphisms $\r^a$ are representable and proper.

Our discussion, in the previous two paragraphs  gives the proof of Theorem 
\ref{thmintro} of the introduction  when $h=1$. 

\begin{Remark} {\rm Of course, the above actually shows that the obvious generalization  of Theorem \ref{thmintro} to $\Cc_{h,d, K}$ for any $h\geq 1$  is also valid.}
\end{Remark}

 \begin{Remarknumb}\label{rem31}
{\rm a) The morphism $\r: \Cc_{h, K}\rightarrow \mathcal R$ is ind-representable: Indeed, suppose that we are given a point $\xi: S=\Spec(R)\rightarrow \RR$ corresponding to a  module $(M, \Phi)$ over $R_W((u))$. Then by Theorem \ref{grass}, the fiber 
$\Cc\times_{\RR,\xi} S\rightarrow S$ is ind-represented
by an ind-algebraic space. We can see that the subspace $\Cc_{h, K}\times_{\RR,\xi} S\hookrightarrow \Cc\times_{\RR,\xi} S$ is described by a closed condition and
it is a closed ind-algebraic subspace. 

b) In general,  
 $\Cc_{h, K}\times_{\RR,\xi} S\rightarrow S$ is not representable
for all $S=\Spec(R)$.
However, assume that $R\simeq \ds{\lim_{\leftarrow {a}}  R/p^aR}$
is a Noetherian $p$-adic ring
and that 
$$
\hat\xi=(\xi^a)\ ,\quad   
\xi^a: \Spec(R/p^aR)\rightarrow \RR^a\ , 
$$
is a point of the formal completion $\widehat\RR$.
Assume also that 
the $R_W\{\{u\}\}$-module $\hat M$ which corresponds to $\hat\xi$ is free over $R_W\{\{u\}\}$. Fix a basis $R_W\{\{u\}\}^d\simeq \hat M$ and set $M=R_W((u))^d\subset \hat M$. The affine Grassmannian $Gr_M\rightarrow S$ is ind-projective
and supports a natural line bundle whose restriction 
on each closed subscheme is very ample. 
As above, we can see that for each $a$, the 
fiber of  the morphism $\r^a:\Cc^a_{h, K}\rightarrow \mathcal R^a$ 
over $\xi^a$ is representable by a closed 
(and hence) projective subscheme of ${\rm Gr}_M\times_{\Z_p}\Z/p^a\Z$. Varying $a$, this defines a formal scheme  over ${\rm Spf}(R)$.
As in \cite{K3}, proof of Prop. 1.3, by using the above ample line bundle on ${\rm Gr}_M$, we may algebraicize this $p$-adic formal scheme over ${\rm Spf}(R)$ to a projective scheme $\Cc_{K,\xi}$ over $\Spec (R)$. The result is that, in this case, the fiber $\widehat\Cc_{d, h, K}\times_{\widehat\RR,\widehat\xi} \hat S \rightarrow \hat S$ between the
formal completions is representable by the formal scheme over $\hat S={\rm Spf}(R)$
associated to the projective scheme $\Cc_{K,\xi}\rightarrow \Spec (R)$. 
If in addition to the hypotheses above, $R/pR$ is of finite type over ${\Bbb F}_p$,
the arguments in the proof of \cite{K3} Prop.\ 1.6.4, show
that the morphism $\Cc_{K,\xi}\rightarrow \Spec (R)$
induces a closed immersion between the generic fibers.

}
\end{Remarknumb}

\subsection{}  \label{cochar}
Choose a cocharacter 
$$
\mu:  {\bar \Q_p}^\times\rightarrow ({\Res}_{ K/\Q_p}GL_d)(\bar \Q_p)=\prod_{\psi: K\rightarrow \bar\Q_p} GL_d(\bar \Q_p)
$$
defined over $\bar \Q_p$ whose conjugacy class is defined over the reflex field $E$. 
The projection to the component corresponding to $\psi$ 
$$
\mu_\psi=pr_\psi\circ \mu: \bar\Q_p^\times\rightarrow GL_d(\bar\Q_p)
$$
provides a grading 
\begin{equation}\label{grading}
\bar\Q_p^d= \bigoplus_{n \in \Z}V^\psi_{n } \ ,
\end{equation}
with $V^\psi_{n  }=\{v\in \bar\Q_p^d\ |\ \mu_\psi(a)=a^{n }v\}$.
Let $h_+$, resp. $h_-$, the maximum, 
resp. minimum value of $n $ (among all the values for all $\psi$) 
for which $V^\psi_{n }\neq (0)$. Set $h=h_+-h_-$. 

We will now  define the corresponding local model 
$M^{\rm loc}_{ \mu, K}$  over $\O_E$; it is going to be a 
projective subscheme of $M_{d,h,K}$ (see \ref{local}).

First, we define the generic fiber of $M^{\rm loc}_{ \mu, K}$
over the reflex field $E$. Suppose that $R$ is a $\bar\Q_p$-algebra, 
and fix an embedding $\psi: K\rightarrow \bar\Q_p$, this induces a homomorphism
$$
R_W=W\otimes_{\Z_p}R\rightarrow R\ ;\quad  a\otimes r\mapsto \psi(a)r\ .
$$
Elements of $M_{h, K}(R)$ correspond 
bijectively to $R_W[u]$-modules $\M$ such that
$$
E(u)^{h_+ } R_W[u]^d\subset  \M\subset E(u)^{h_- } R_W[u]^d
$$
with graded pieces $R_W$-projective and with rank locally constant on $\Spec(R)$. Write
$$
{\rm Norm}_{K_0/\Q_p}(E(u))=\prod_{\psi: K\rightarrow\bar\Q_p}(u-\varpi_\psi)\in \bar\Q_p[u]
$$
so that $\varpi_\psi=\psi(\pi)$. 
Using this, we see that we can write $\M=\oplus_\psi \M_\psi$
with $\M_\psi$  a $R[u]$-submodule with 
$$
(u-\varpi_\psi)^{h_+}R[u]^d\subset \M_\psi\subset (u-\varpi_\psi)^{h_-}R[u]^d.
$$
For each such $\psi$,  consider the $R$-module  
$$
\M_\psi\cap (u-\varpi_\psi)^j R[u]^d/
\M_\psi\cap (u-\varpi_\psi)^{j+1}R[u]^d
$$
We ask that for each $\psi$,  
$j\in \Z$, this is a projective $R$-module of rank $\dim(V^\psi_j)$.
 We can see that this condition defines a locally closed subvariety 
of $M_{h, K}\otimes_{\Z_p}\bar\Q_p$.
This carries an action of ${\rm Gal}(\bar\Q_p/E)$ 
that allows us to descend it to a   subvariety $Z$
of $M_{h, K}\otimes_{\Z_p}E$. By definition, the generic fiber of the 
local model $M^{\rm loc}_{\mu, K}$ is the Zariski closure $\bar Z$ of $Z$
in $M_{h, K}\otimes_{\Z_p}E$. Finally, by definition,  
the local model $M^{\rm loc}_{\mu, K}$ is the flat 
closure of $\bar Z$ in $M_{h, K}\otimes_{\Z_p}\O_E$.
\smallskip

Observe that, by the above, for each $\O_E$-scheme $S=\Spec(R)$ in ${\rm Nil}_p$ 
  we have
\begin{equation}
 M^{\rm loc}_{\mu, K}(R)\hookrightarrow (\F_G\otimes_{\Z_p}\O_E)(R).
\end{equation}

\begin{Remarknumb}\label{mini}
{\rm Suppose that, for all $\psi$, we   have $n\in \{0,1\}$ 
in (\ref{grading}). Then $\mu$ is miniscule.  Assume $h_+=1$, $h_-=0$, which is the typical case.
Then $h=1$,  $R_W[u]/(E(u)^h)=\O_K\otimes_{\Z_p}R $,
and $M_{h,K}(R)$ is given by $\O_K\otimes_{\Z_p}R$-submodules
$$
\E\subset (\O_K\otimes_{\Z_p}R)^d
$$
which are locally on $R$ direct summands as $R_W$-modules.
Set $r_\psi=\dim(V^\psi_0)$.
The conditions above amount to asking that
${\rm rank}_R(\E_\psi)=r_\psi$ and so $M^{\rm loc}_{\mu, K}$ agrees with the local model of [PR1]. 
If $k=\Ff_p$, the special fiber $M^{\rm loc}_{\mu, K}\otimes_{\O_E}\bar\Ff_p$ can be identified with the affine Schubert variety 
$S_\nu$ in the affine Grassmannian of $GL_d$, where the coweight $\nu$ is the dual partition to $(r_\psi)_\psi$, i.e.,
$
\nu_i=\sharp \{\psi\mid r_\psi\geq i\}\ ,
$
cf. \cite{P-R1}, Thm.\ 5.4. }
\end{Remarknumb}

\subsection{}\label{3c}
We continue with the above notations of \S \ref{local}.   By definition we have a closed immersion
$$
\widehat M^{\rm loc}_{\mu, K}\hookrightarrow \widehat M_{h, K}
$$
which is equivariant for the natural action of (the formal completion of) $L^+G$. 
Using descent and (\ref{qhK}) we obtain a fpqc stack $\widehat\Cc_{\mu, K}$
over ${\rm Nil}_p\cap ({\rm Sch}/\O_E)$ together with a formally smooth morphism
$$
q_{\mu, K}: \widehat\Cc_{\mu, K}\rightarrow [L^+G\backslash \hat M^{\rm loc}_{\mu, K}].
$$
Indeed, if $S=\Spec(R)$ is an $\Spec(\O_E)$-scheme in ${\rm Nil}_p$ then 
$\widehat\Cc_{\mu, K}(R)$ is the groupoid of pairs
$(\Mfr, \Phi)$
of a $R_W[[u]]$-module $\Mfr$ which is   
$R_W[[u]]$-free of rank $d$ (fpqc) locally on $R$ and a 
$R_W((u))$-module isomorphism 
\begin{equation}
\Phi : \phi^*{\frak M}[1/u] \xrightarrow{\sim} {\Mfr}[1/u]\ , 
\end{equation} 
such that, locally, there is an isomorphism $\alpha: R_W[[u]]^d\xrightarrow{\sim} \Mfr$
for which the $R_W[[u]]$-lattice $\alpha^{-1}(\Phi(\phi^*\Mfr))\subset R_W((u))^d$
belongs to the subset $\hat M^{\rm loc}_{\mu, K}(R)$ of  $(\F_G\otimes_{\Z_p}\O_E)(R)$.
As in \S \ref{3b2}, we see that there is also a morphism of formal completions
\begin{equation}
 \widehat\r_{\mu}: \widehat\Cc_{\mu, K}\rightarrow \widehat{\RR}_d\otimes_{\Z_p} \O_E\ 
\end{equation}
which can be obtained as the limit of representable and proper
morphisms.

\section{Deformations of Galois representations}\label{3}

In this section, we explain an aspect of the connection between 
the spaces of $\Phi$-modules and the deformation theory of Galois representations
as developed by Kisin \cite{K1}, \cite{K3}. We restrict attention
to the {\sl flat} or {\sl Barsotti-Tate} case (cf. \cite{K1})
This corresponds to the case $h=1$. For simplicity, we also assume $p$ is odd.

\subsection{Galois representations} 
 
Suppose that $R=\Lambda$ is a  $\Z_p$-algebra with finitely many elements.
As in \S 1 (see also \cite{Fontaine}), a pair $(M, \Phi )$ corresponding to an object of ${\cal R}(\Lambda)$ gives
 an \'etale $\Lambda$-sheaf over $\Spec(k((u)))$ which is free of rank $d$,
i.e an equivalence class of a representation 
$$
\rho_{(M, \Phi)}: {\rm Gal}(k((u))^{\rm sep}/k((u)))\rightarrow GL_d(\Lambda).
$$
As a result, an object $(\Mfr, \Phi )$ of 
${\cal C}_d(\Lambda)$ also gives a representation
$\rho_{(\Mfr[1/u], \Phi)}$ of the Galois group ${\rm Gal}(k((u))^{\rm sep}/k((u)))$.

Now let $\Ff$ be a finite field and suppose that
$$
\rho: {\rm Gal}(k((u))^{\rm sep}/k((u)))\rightarrow GL_d({\Ff})
$$
is a representation which corresponds to a pair $(M_0,\Phi_0)$, and 
 consider the corresponding object $[\rho]: \Spec(\Ff)\rightarrow {\cal R}$.
Denote by ${\cal R}_{[\rho]}$ the groupoid over finite local Artinian $\Z_p$-algebras $\Lambda$ 
with residue field $\Ff$, with objects
$$
{\cal R}_{[\rho]}(\Lambda)=\{(M,\Phi)\in {\cal R}(\Lambda), \alpha: (M,\Phi)\otimes_\Lambda {\Ff}\xrightarrow{\sim}  (M_0,\Phi_0)\}
$$
and obvious morphisms. By the above, ${\cal R}_{[\rho]}$ is identified with the groupoid of deformations ${\frak D}_\rho$ of the Galois representation $\rho$.

\subsection{Finite flat group schemes}

 Now let $K$ be a finite extension of $\Q_p$ with residue field $k$ and ramification index $e$.
Choose a uniformizer $\pi$ of $K$ with Eisenstein polynomial $E(u)$ 
over ${\rm Fr}(W)$. Set $K_\infty=\cup_nK(\pi_n)$, where $\pi_n=\pi^{1/p^n}$
are compatible choices of roots; then the theory of norm fields allows us to identify the Galois groups
$$
G_\infty:={\rm Gal}(\bar K/K_\infty)\xrightarrow{\sim} {\rm Gal}(k((u))^{\rm sep}/k((u)))\ ,
$$
comp. \cite{K1}, \S 1. The following can be derived from \cite{K2} Theorem 0.5
by taking into account the functoriality of the $\Lambda$-action 
and the properties of the Breuil-Kisin module functors (for example see \cite{K1} \S 1.2).

\begin{thm} (Kisin)
Assume $p>2$ and let $\Lambda$ be a ${\Z}_p$-algebra with finitely many elements.  There is 
an equivalence   between the groupoid of finite flat commutative group schemes $\cal G$ with an action of $\Lambda$ (i.e ``$\Lambda$-module schemes")
over $\O_K$ such that $\Gg(\bar\O_K)\simeq \Lambda^d$, and the groupoid of pairs   $(\Mfr, \Phi)$
of $\Lambda_W[[u]]$-$\Phi$-modules 
 with
the following properties:

a)  ${\rm Coker}(\Phi)$ is annihilated by $E(u)$,

b)  $\Mfr[1/u]$ is $\Lambda_W((u))$-free of rank $d$,

c) $\Mfr$ is a $W[[u]]$-module
of projective dimension $1$, i.e., equivalently by \cite{K2}, Lemma (2.3.2),  
$\Mfr$ is an iterated extension of free $k[[u]]$-modules.

\smallskip

Under this equivalence,  the restriction of 
$$
\rho_{\Gg}: {\rm Gal}(\bar K/K)\rightarrow {\rm Aut}_{\Lambda}(\Gg(\bar\O_K))
$$
  to $G_\infty\simeq {\rm Gal}(k((u))^{\rm sep}/k((u)))$
is isomorphic to $\rho_{(\Mfr[1/u],\Phi)}(1)$\ {\rm [twist by the cyclotomic character]}. 
\end{thm}

Note that property (c)  is automatically satisfied  when $p\cdot \Lambda=(0)$. We will also consider the groupoid of modules as above with 
the additional property 
\smallskip

{\it d) $\Mfr$ is $\Lambda_W[[u]]$-free}.
\smallskip

\begin{lemma} \label{41}
Let $A$ be a local Noetherian ring with residue field $l$ and suppose that 
$M\subset  A((u))^d$ is a finitely generated $A[[u]]$-module such that $M[1/u]= A((u))^d$ and $M\otimes_{A}l\simeq l[[u]]^d$. Then $M\simeq A[[u]]^d$.
\end{lemma}

\begin{proof}
By Nakayama's lemma, there is a surjective $A[[u]]$-homomorphism $\phi: A[[u]]^d\rightarrow M$. 
Then $\phi[1/u]: A((u))^d\rightarrow M[1/u]=A((u))^d$ is also a surjection
which then has to be bijective. This implies that $\phi$ is also 
injective.
\end{proof}
\smallskip

Write $\Lambda\otimes _{\Z_p}W=\prod_j \Lambda_j$ with $\Lambda_j$ Artin local. An application of the lemma to $\Lambda_j$ shows that
the additional property $(d)$ is satisfied if and only if $\Mfr\otimes_{\Lambda}{\Ff}\simeq ({\Ff}\otimes k)[[u]]^d$.

\subsection{} 
Suppose that $(A, {\mathfrak m})$ is an Artin local Noetherian ring 
with finite residue field $\Ff$. A representation  
$$
\rho: {\rm Gal}(\bar K/K)\rightarrow GL_d(A)
$$
is called {\sl flat} \cite {Ra} if the corresponding
$\Z_p[{\rm Gal}(\bar K/K)]$-module  
is isomorphic to the twist by $(-1)$ of the module obtained by the Galois action on the ${\Z}_p$-module of $\O_{\bar K}$-points of some 
commutative finite flat group scheme over $\O_K$.
This notion extends to the more general situation that $(A, {\mathfrak m})$ is a  complete local Noetherian ring 
with finite residue field $\Ff$. In this case, the representation 
$$
\rho: {\rm Gal}(\bar K/K)\rightarrow GL_d(A)
$$
is {\sl flat} iff, for all $n\geq 1$, the representation obtained 
by reducing $\rho$ modulo ${\mathfrak m}^n$ is flat, cf. \cite{Ra}.
%A representation of $G_\infty={\rm Gal}(k((u))^{\rm sep}/k((u)))$
%is $K$-flat if it is the restriction of a flat representation of ${\rm Gal}(\bar K/K)$ under $G_\infty\hookrightarrow {\rm Gal}(\bar K/K)$.

Consider the   morphism of formal stacks 
$
\r_K: \widehat\Cc_{K}:=\widehat\Cc_{1,K}\rightarrow \widehat\RR
$. Suppose that $(A, {\mathfrak m})$ is a  complete local Noetherian ring 
with finite residue field $\Ff$. 
Let $\xi=(\xi_n)_{n\geq 1}\in \widehat\RR(A)$ be 
an $A$-valued object of $\widehat\RR$, 
where for each $n\geq 1$, $\xi_n$ is in $\widehat\RR(A/{\frak m}^n)$.
The $2$-fiber product $\xi_n\times _{\widehat\RR}\widehat\Cc_K$ 
is representable by a  projective scheme
$
\Cc_{K,\xi_n}\rightarrow \Spec(A/{\frak m}^n)
$.
In the limit, we obtain a formal scheme $\Cc_{K, \xi}$
over ${\rm Spf}(A)$. The argument in \cite{K3}, Cor. 1.5.1 (or see Remark \ref{rem31} (b)) shows that 
this is algebraizable to a projective scheme
 $$
\Cc_{K,\xi }\rightarrow \Spec(A ).
$$
Denote by $A^K$ the quotient of $A$ that 
corresponds to the scheme theoretic image of this morphism.
We obtain 
$$
\xi_K: \Spec(A^K)\rightarrow \Spec(A)\rightarrow \widehat\RR.
$$
\begin{prop}
 With the above notations, 
assume in addition that $A$ is Artinian. Then 
$$
\rho_{\xi_K}: G_\infty={\rm Gal}(k((u))^{\rm sep}/k((u)))\rightarrow
GL_d(A^K)
$$
extends to a representation of ${\rm Gal}(\bar K/K)$
which is flat.
\end{prop}
 
\begin{proof}
For simplicity, set   $B=A^K$. Then there is $B\hookrightarrow B'$ with $B'$ a $B$-algebra of finite type
such that $\Cc_{K}$ affords a $B'$-valued point $\zeta'$ that lifts $\xi_B$.
Denote by $(M_B, \Phi)$ the $B_W((u))$-module that corresponds to $\xi_B$.
Since $B$ is Artinian, $M_B\simeq B_W((u))^d$. 
Giving the point $\zeta'$ amounts to giving a $B'_W[[u]]$-projective module $\Mfr'$
of rank $d$ in $M'=M_B\hat\otimes_BB'$
which satisfies
 $$
E(u)\Mfr'\subset \Phi(\phi^*\Mfr')\subset \Mfr'.
 $$
Now set $\Mfr:=\Mfr'\cap M\subset M'$; this is a $B_W[[u]]$-module which is $\Phi$-stable. 
The proof of [Ki1] Prop. 2.1.4 applies to show that $\Mfr$ satisfies properties (a), (b) and (c)
(with $\Lambda=B$). Therefore, the $\Z_p[G_\infty]$-module given by $\rho(1)$ is given by
a finite flat group scheme over $\O_K$ as desired. (However, note here that, as pointed out by Kisin, $\Mfr$ does not have to be
$B_W[[u]]$-free.)
\end{proof}

\subsubsection{}  Assume now that $(A, {\mathfrak m})$ is a  complete local Noetherian ring 
with finite residue field $\Ff$ and that in addition:

1) the $A$-valued  representation $\rho_\xi$ of
$G_\infty$ which corresponds to
$\xi$ extends to
a representation of ${\rm Gal}(\bar K/K)$, and

2) the $\Ff$-valued  representation of
${\rm Gal}(\bar K/K)$ which is obtained by reducing 
$\rho_\xi$ modulo $\mathfrak m$ is flat.

By [Ra] there is a quotient $A\rightarrow A^{\rm fl}$
such that for each $A\rightarrow B$ with $B$ a local Artinian $A$-algebra with residue field $\Ff$, 
the representation obtained by composing with $GL_d(A)\to GL_d(B)$ is flat if and only if $A\rightarrow B$ 
factors as $A\rightarrow A^{\rm fl}\rightarrow B$. 
Using the above proposition one can see that $A\rightarrow A^K$ factors as
$$
A\rightarrow A^{\rm fl}\rightarrow A^K.
$$

\begin{Remarknumb}\label{defgal}
{\rm In general, it is not clear  whether we should expect that
$A^{\rm fl}\rightarrow A^K$ is an isomorphism.  
The issue is the following: Consider a deformation $\rho$ of $\rho_0=$
$\rho_\xi$ modulo $\mathfrak m$
over a finite Artin local $\Z_p$-algebra $\Lambda$ with residue field $\Ff$.
It corresponds to
a $\Lambda_W((u))$-$\Phi$-module $M$. By assumption,  the $(\Ff\otimes_{\Z_p} W)((u))$-$\Phi$-module
$M_0=M\otimes_\Lambda \Ff$ which corresponds to $\rho_0$ contains a $(\Ff\otimes_{\Z_p} W)[[u]]$-$\Phi$-submodule
$\Mfr_0$ with $E(u)\Mfr_0\subset \Phi(\phi^*\Mfr_0)\subset \Mfr_0$. We can easily see 
that $\Mfr_0\simeq (\Ff\otimes_{\Z_p} W)[[u]]^d$. Assume now that $\rho$ is also flat;
this implies that $M$ contains a $\Lambda_W[[u]]$-$\Phi$-module $\Mfr$ with $\Mfr[1/u]=M$
that satisfies properties  (a), (b) and (c). The problem is
that if $e>p-1$, we cannot expect that $\Mfr$ is a deformation of $\Mfr_0$
(so that we can apply Lemma \ref{41}).
The question is: Is there is some $\Z_p$-algebra $C$ that contains $\Lambda$ and a $C_W[[u]]$-$\Phi$-module $\Mfr_C$ 
with $\Mfr_C[1/u]=M\otimes_\Lambda C$ which is in addition 
$C_W[[u]]$-projective and such that
$\Mfr=M\cap \Mfr_C$?
Our discussion implies}
\end{Remarknumb}
\begin{prop}
In the situation of Remark \ref{defgal}, assume that $e\leq p-1$. Then
$
A^{\rm fl}\simeq A^K\ .
$
\end{prop}

\medskip

\section{Coefficient domains and a period morphism}\label{dom}

\subsection{}
Fix $d$, the local field $K$ and $h\geq 1$.
We define the stack in groupoids ${\mathcal D}_{d, h, K}$ over schemes
over $\Q_p$ which is described as follows:

If $R$ is a $\Q_p$-algebra, then the objects of ${\mathcal D}_{d,h, K}
(R)$ are triples $(D, \Phi, {\rm Fil}^\bullet)$
where

   \begin{itemize}
\item $D$ is a $R\otimes_{\Q_p} K_0$-module which is, locally on $R$,
free of rank $d$,
\item $\Phi: D\rightarrow D$ is an ${\rm Id}\otimes_{\Q_p} \phi$-linear
automorphism,
\item ${\rm Fil}^\bullet$ is a exhausting, decreasing filtration of
$D_K:=D\otimes_{K_0}K$
by $R\otimes_{\Q_p}K$-modules which are locally direct summands
and satisfy ${\rm Fil}^0=D_K$, ${\rm Fil}^{h+1}=(0)$.
\end{itemize}

We can see that  ${\mathcal D}_{d,h, K}$ is a fpqc stack over $\Q_p$.
Locally we can choose a basis of $D$;
this allows us to write the stack as a quotient
\begin{equation}
{\mathcal D}_{d, h, K}=[({\rm Res}_{K_0/\Q_p}GL_d\times_{\Q_p} {\rm Gr}_{d, h, K})/_{(\phi, \cdot )} {\rm Res}_{K_0/\Q_p}GL_d ].
\end{equation}
Here ${\rm Gr}_{d, h, K}$ is the Grassmannian of filtrations as above
of length $h+1$ on a vector space of dimension $d$ over $K$.
(Notice here that we are not yet prescribing  dimensions 
for the graded pieces ${\rm Fil}^i/{\rm Fil}^{i+1}$; in particular,
${\rm Gr}_{d, h, K}$ is not necessarily connected.)
The symbol $(\phi, \cdot)$ is supposed to
remind us that
the action of the group ${\rm Res}_{K_0/\Q_p}GL_d $ on
the product is by $\phi$-conjugation on the first factor and by
translation on the second.
It follows from this description that ${\mathcal D}_{d, h, K}$ is an Artin stack, smooth of
finite type over $\Q_p$. 

Similarly, suppose that  $\mu:{\bar \Q_p}^\times\rightarrow ({\Res}_{ K/\Q_p}GL_d)
(\bar \Q_p)$ is a coweight as
in \S \ref{3c} before. Assume that $\mu$ is defined over the
reflex field $E$ and, for simplicity, assume
$h_-=0$ so that $h=h_+$. We define the stack in groupoids ${\mathcal D}_{\mu, K}$ over schemes
over $E$ which is described as follows:
If $R$ is an $E$-algebra, then the objects of ${\mathcal D}_{\mu, K}
(R)$ are triples $(D, \Phi, {\rm Fil}^\bullet)$ that correspond to objects of ${\mathcal D}_{d,h, K}(R)$
 as above with the additional property 

   \begin{itemize}
\item The filtration 
${\rm Fil}^\bullet$ is
of type $\mu$ in the sense that the base change
of the graded piece $ {\rm Fil}^j/{\rm Fil}^{j+1}$ under ${\rm id}\otimes\psi: R\otimes_{\Q_p}K\rightarrow
R\otimes_{\Q_p}\bar\Q_p$ has rank equal to $\dim(V_j^\psi)$ for each $\psi$ and each $j$ (see (\ref{grading})).
\end{itemize}

Once again ${\mathcal D}_{\mu, K}$ is a fpqc stack over $E$
which is an Artin stack, smooth of
finite type over $E$. 
We can write
\begin{equation}
{\mathcal D}_{\mu, K}=[({\rm Res}_{K_0/\Q_p}GL_d)_E\times_E {\rm Gr}_
{\mu, K}/_{(\phi, \cdot )}({\rm Res}_{K_0/\Q_p}GL_d)_E].
\end{equation}
Here ${\rm Gr}_{\mu, K}$ is the Grassmannian of filtrations as above
of type $\mu$
on $(E\otimes_{\Q_p}K)^d$.

\begin{Remarknumb}
{\rm We can also consider the following ``rigid" variants:   ${\mathfrak D}_{d, h, K}={\mathcal D}_{d, h, K}^{\rm rig}$
is the category fibered in groupoids over the category of rigid spaces over $\Q_p$ which is defined as follows. If ${\mathfrak X}$ is a rigid space, then ${\mathfrak D}_{d, h, K}({\mathfrak X})$ is the groupoid of pairs $(D, \Phi, {\rm Fil}^\bullet)$ with
$D$ a coherent sheaf of $\O_{\mathfrak X}\otimes_{\Q_p} K_0$-modules over $\mathfrak X$ which is locally  
free of rank $d$, $\Phi$ an $1\otimes_{\Q_p}\phi$-linear isomorphism of $D$,
and ${\rm Fil}^\bullet$ a filtration of $D\otimes_{K_0}K$ (of length $h$, as above) by 
coherent $\O_{\mathfrak X}\otimes_{\Q_p} K$-sheaves   over ${\mathfrak X}$
with locally free graded pieces.  Here we are implicitly using that descent  of coherent modules under fpqc morphisms of rigid spaces is effective, cf. \cite{BG}. Similarly, we can define ${\mathfrak D}_{\mu, K}$
etc. 
}\end{Remarknumb}

\subsubsection{} \label{5a1}
Now suppose that $\A$ is a $p$-adic ring.
Set $A=\A[1/p]$. 
Let $(\Mfr, \Phi)$ a $\A_W[[u]]$-$\Phi$-module which corresponds
to an object
of $\widehat\Cc_{d,h, K}(\A)$. Consider
$$
D=(\Mfr/u\Mfr)[1/p]
$$
with the $1\otimes \phi$-endomorphism given by $\Phi\ {\rm mod\ } u$;
we can easily see that this endomorphism is bijective and that
$D$ is an $A\otimes_{\Q_p}K_0$-module which is projective  of rank $d$.

Similarly to [Ki2], [Ki3], we set
  $$
  \O_A:=\lim_{\leftarrow n}  (\A_W[[u, u^n/p]] [1/p])\subset
A_W[[u]]=(A\otimes_{\Q_p}K_0)[[u]];
  $$
  in particular
  $$
  \O=\O_{\Q_p}:=\lim_{\leftarrow n}  (W[[u, u^n/p]][1/p])\subset
K_0[[u]]
  $$
  is the ring of rigid analytic functions on the open disk $\Bbb U$ of radius
$1$ over $K_0$.
  (The inverse limit is under the maps given by $u^{n'}/p\mapsto u^{n'-
n}\cdot (u^n/p)$.)

  We have inclusions $W[[u]][1/p]\hookrightarrow \O$, $\A_W[[u]]
[1/p]\hookrightarrow \O_{A}$.
  The endomorphism $\phi$ has a unique continuous extension to
  $\O$ and $\O_{A}$.  Set $\M=\Mfr\otimes_{\A_W[[u]]}\O_{A}$;
the $\Phi$-structure on $\Mfr$ induces
$$
\Phi: \phi^*(\M)\rightarrow \M.
$$
This map is injective and we have $E(u)^h\M\subset \Phi(\phi^*(\M))
\subset  \M$, and $D=\M/u\M$.
  Set
  \begin{equation*}
  \lambda=\prod_{n=0}^\infty \phi^n(E(u)/E(0))\in \O.
  \end{equation*}
For each $m$, let $r(m)$ be the smallest integer such that
  $em<p^{r(m)}$ and consider
  $$
\O_{\A, e, m}:=\A_W[[u, u^{p^{r(m)}}/p^m]] [1/p].
  $$
  There is a  ring homomorphism
  $
  \O_{A}\rightarrow \O_{\A, e, m}
  $.
  Since $|\pi|=p^{-1/e}$, we can see that   $u\mapsto \pi$
  gives $\O_{\A, e, m}\rightarrow A\otimes_{\Q_p}K$ which induces an
isomorphism
  \begin{equation}\label{425}
  \O_{\A, e, m}/(E(u)) \xrightarrow{\sim} A\otimes_{\Q_p}K.
  \end{equation}

  Recall $\M=\Mfr\otimes_{\A_W[[u]]}\O_A$. As in [Ki3] Lemma (2.2) we   see that
there is
  a unique $\phi$-compatible $A_W$-linear map
  \begin{equation*}
  \xi: D\rightarrow \M
  \end{equation*}
  with the following properties:
  
1) The reduction modulo $u$ of $\xi$ is the identity. 

2) The induced map $\xi: D\otimes_{A\otimes_{\Q_p}K_0} \O_A
\rightarrow \M$
is injective and has cokernel killed by $\lambda^h$.

3) For any sufficiently large $m$, the induced map
$$
\xi\otimes_{\O_A}  \O_{\A, e,m}:
D\otimes_{A\otimes_{\Q_p}K_0}  \O_{\A, e, m} \rightarrow \M
\otimes_{\O_A}  \O_{\A, e,m}=
\Mfr\otimes_{\A_W[[u]]}  \O_{\A, e, m}
$$
is injective and its image is equal to that
of the map
$$
\phi^*(\M)\otimes_{\O_A} \O_{\A, e, m}\rightarrow \M\otimes_{\O_A}
\O_{\A, e, m}
$$
induced by $\Phi$.

Indeed, the construction of $\xi: D\rightarrow \M$ in [Ki3] works
verbatim when $\Mfr$ is $\A_W[[u]]$-free.   (The assumptions
that $\A$ is complete, local and Noetherian are not needed
for our version of the construction.
We can choose $\Mfr$ and $\Mfr/u\Mfr$ to play the roles of the modules
denoted by $\Mfr_A^\circ$ and $D^\circ_A$ in loc. cit.) The general case
 is obtained by gluing, using the uniqueness of $\xi$ in the free case.
To check claims (2) and (3) we  can argue as in [Ki3].
(Note that loc. cit. Lemma (2.2.1) is also valid with essentially the
same proof even when $\A$ is
not  Noetherian.)

As a result of (2) and (3), and by using (\ref{425}) to reduce modulo $(E(u))$, we obtain
an isomorphism
\begin{equation}\label{isoFil2}
D\otimes_{K_0}K\xrightarrow{\sim} \phi^*(\M)\otimes_{\O_A} (A\otimes_{\Q_p}K)\simeq \Phi(\phi^*\Mfr[1/p])\otimes_{\A_W[[u]][1/p]}  (A\otimes_{\Q_p}K).
\end{equation}
(In the last tensor product, we use $\A_W[[u]][1/p]\rightarrow \A_W[[u]][1/p]/(E(u))\xrightarrow{\sim} A\otimes_{\Q_p}K$.)
We can see that the isomorphism (\ref{isoFil2}) is independent of the choice of $m$.

\subsubsection{} \label{rigid}

Recall that $\Bbb U$ denotes the rigid open unit disk over $K_0$.
If $I$ is a subinterval of $[0,1)$, we set
  ${\Bbb U}(I)$ for the admissible open subspace of points with
  absolute value in $I$. Set $\O_I=\Gamma({\Bbb U}(I), \O_{{\Bbb U}(I)})$
  so that $\O=\O_{[0,1)}$.
We denote by $\phi: {\Bbb U}\rightarrow {\Bbb U}$ the ``Frobenius" morphism
which corresponds to $\phi: \O\rightarrow \O$ as before.

We can  consider the category $\Ccc_{d,h, K} $ fibered in groupoids  over  $\Q_p$-rigid spaces which is defined as follows. 
Let $\mathfrak X$ be a rigid space over $\Q_p$ and consider ${\mathfrak X}\times {\Bbb U}$ with the partial Frobenius $\phi:={\rm id}\times\phi: {\mathfrak X}\times {\Bbb U}\rightarrow {\mathfrak X}\times {\Bbb U}$. Then, by definition,  $\Ccc_{d, h, K}({\mathfrak X})$ is 
the groupoid of pairs $(\M, \Phi)$ where $\M$ is a coherent sheaf over ${\mathfrak X}\times {\Bbb U}$ which is locally on $\mathfrak X$ free of rank $d$, and $\Phi: \phi^*\M\rightarrow \M$ an injective homomorphism
with cokernel annihilated by $E(u)^h$. 

Denote by 
$i: {\mathfrak X}_{K_0} \hookrightarrow {\mathfrak X}\times {\Bbb U}$ the inclusion
$i(x)=(x,0)$ and by ${ p}: {\mathfrak X}\times {\Bbb U}\rightarrow {\mathfrak X}_{K_0}$ the projection. If $(\M, \Phi)$ is an object of $\Ccc_{d, K} ({\mathfrak X})$ we set $D=i^*\M$. This is a coherent sheaf on $\mathfrak X_{K_0}$ which is locally free of rank $d$; the morphism $\Phi: \phi^*\M\rightarrow \M$ induces a $\phi$-linear isomorphism $\Phi:  D\rightarrow D$.

\begin{prop}\label{xirigid}
There is a (unique) $\Phi$-compatible morphism of sheaves of $\O_{{\mathfrak X}_{K_0}} $-modules
$
\xi: D\rightarrow {p}_*(\M)
$
such that

1) $i^*\xi$ is the identity,

2) the induced morphism ${ p}^*\xi: { p}^*D\rightarrow 
 \M$ is 
injective and has cokernel annihilated by $\lambda^h$,

3) If $r\in (|\pi|, |\pi|^{1/p})$, then the image of 
the restriction ${ p}^*\xi_{[0,r)}$ to ${\mathfrak X}\times {\Bbb U}[0,r)$ coincides with 
the image of $\Phi_{[0,r)}: \phi^*\M_{[0,r)}\rightarrow \M_{[0,r)}$.
\end{prop}

\begin{proof}
When ${\mathfrak X}={\rm Sp}(\Q_p)$ is a point, this is \cite{K2} Lemma 1.2.6.
Note that there is at most one $\Phi$-compatible $\xi: D\rightarrow p^*(\M)$ 
that satisfies  property (1).
Indeed, if $\xi$, $\xi'$ are two such  morphisms we have ${\rm Im}(\xi-\xi')\subset u\cdot p^*(\M)$.
The $\Phi$-compatibility gives $\Phi\cdot (\xi-\xi')=(\xi-\xi')\cdot \Phi$.
Hence, since $\Phi: D\rightarrow D$ is an isomorphism, we obtain inductively
${\rm Im}(\xi-\xi')\subset u^s\cdot p^*(\M)$ for all $s\geq 0$. This implies that $\xi=\xi'$.
To show the existence of $\xi$ we
suppose first that $A$ is a Tate $\Q_p$-algebra
and that ${\mathfrak X}={\rm Sp}(A)$ is the corresponding affinoid rigid space.
Then $\O_A$ is the ring of rigid analytic functions $\O_{{\rm Sp}(A)\times {\Bbb U}}$ on the product ${\rm Sp}(A )\times {\Bbb U}$ and $\M$ is given by an $\O_A$-module as in the previous paragraph. 
There is a $p$-adic ring $\mathcal A$ which is topologically of finite presentation (tfp) over $\Z_p$ and $p$-torsion 
free such that $A=\A[1/p]$. The arguments of \cite{K2} Lemma 1.2.6,
\cite {K3} Lemma (2.2) (see also the previous paragraph) extend to this case to construct $\xi$
that satisfies all the required properties.
The result in the case of a general rigid space $\mathfrak X$ follows by the affinoid case
above by gluing using the uniqueness of $\xi$.
\end{proof}
\smallskip

Suppose now that $r$ is in $(|\pi|, |\pi|^{1/p})$. Then $\O_{{\mathfrak X}\times {\Bbb U}[0,r)}/(E(u))\simeq \O_{\mathfrak X}\otimes_{\Q_p}K$.
As a result of (2) and (3), we obtain an isomorphism
\begin{equation}\label{isoFil3}
p^*D_{[0,r)}\xrightarrow{\ \sim\ }
\Phi(\phi^*\M)_{[0,r)}\
\end{equation}
which by reducing modulo $(E(u))$ gives an isomorphism 
\begin{equation}\label{isoFil4}
D\otimes_{K_0}K\xrightarrow{\ \sim\ }
\Phi(\phi^*\M)\otimes_{\O_{{\mathfrak X}\times {\Bbb U}[0,r)}}(\O_{\mathfrak X}\otimes_{\Q_p}K)
\end{equation}
of coherent $\O_{\mathfrak X}\otimes_{\Q_p}K$-sheaves over $\mathfrak X$.

\subsection{} In what follows, we will assume that  $h=1$.

Let $\A$ be a $p$-adic ring
and set $A=\A[1/p]$. 
Suppose that $(\Mfr, \Phi)$ is 
 an $\A_W[[u]]$-$\Phi$-module  which corresponds
to an object
of $\widehat\Cc_{d, K}(\A)$. We will define
an object $D(\Mfr, \Phi)$ of ${\mathcal D}_{d, K}(A)$ by
following the construction
of [Ki2], [Ki3]. 
Let
$
D=(\Mfr/u\Mfr)[1/p]
$
with its $\Phi$-structure be as in \S \ref{5a1}. It  remains to define the filtration ${\rm Fil}^\bullet$ on $D
\otimes_{K_0}K$.
Since $h=1$, we have
$$
E(u)\Mfr\subset \Phi(\phi^*\Mfr)
\subset  \Mfr\ .
$$
 The module $\Phi(\phi^*\Mfr[1/p])$ 
is filtered
$$
 E(u) \Phi(\phi^*\Mfr[1/p]) \subset E(u)\Mfr[1/p]\subset \Phi(\phi^*\Mfr[1/p])\ .
 $$
Hence, we can filter the $A\otimes_{\Q_p}K$-module
  $D\otimes_{K_0}K\simeq \Phi(\phi^*\Mfr[1/p])\otimes_{\A_W[[u]][1/p]}  (A
\otimes_{\Q_p}K)$ via (\ref{isoFil2}) by  taking the image of this filtration, i.e we set
\begin{eqnarray*}
{\rm Fil}^2&=&(0),\\
{\rm Fil}^1&= & E(u)\Mfr[1/p]{\,} {\rm mod}\, E(u)\Phi(\phi^*\Mfr[1/p]),\\
{\rm Fil}^0&=& D\otimes_{K_0}K\ .
\end{eqnarray*}
Since $E(u)$ is not a zero divisor in $\A_W[[u]][1/p]$, we can see
(cf. \cite{K3} 2.6.1 (1)) that the quotient $\Mfr[1/p]/\Phi(\phi^*\Mfr[1/p])$ is a finitely generated
projective  $A\otimes_{\Q_p}K$-module. We conclude that
 ${\rm Fil}^1$  is a finitely generated projective $A\otimes_{\Q_p}K$-module
 which is locally a direct summand of $D\otimes_{K_0}K$. Hence, 
 $(D, \Phi, {\rm Fil}^\bullet)$ gives an object of ${\mathcal D}_{d, K}(A)$
 which we will denote by $ D(\Mfr, \Phi)$. This gives a functor of groupoids, 
 $$
 \underline D(\A): \widehat\Cc_{d, K}(\A)\rightarrow  {\mathcal D}_{d, K}(A).
 $$

 \subsubsection{} Similarly, if $(\M,\Phi)$ is an object of $\Ccc_{d,K}({\mathfrak X})$
 for a rigid space $\mathfrak X$, we  consider  $D=i^*\M$ with its $\phi$-linear isomorphism $\Phi:  D\rightarrow D$
 as above. We also have
 $$
 E(u) \Phi(\phi^*\M) \subset E(u)\M\subset \Phi(\phi^*\M)\ 
 $$
(A filtration of coherent sheaves
over ${\mathfrak X}\times {\Bbb U}$.) As above, we can use this and (\ref{isoFil4})
to produce a filtration 
$$
(0)={\rm Fil}^2\subset {\rm Fil}^1\subset {\rm Fil}^0=D\otimes_{K_0}K
$$
of the coherent sheaf $D\otimes_{K_0}K$ over  ${\mathfrak X}_{K_0}$.
We can see that the triple $(D, \Phi, {\rm Fil}^\bullet)$ gives an object of ${\mathfrak D}_{d, K} ({\mathfrak X})$. Since  $(\M, \Phi)\mapsto (D, \Phi, {\rm Fil}^\bullet)$ is functorial, this defines a functor of groupoids, ${\Ccc}_{d, K} ({\mathfrak X}) \rightarrow  {\mathfrak D}_{d, K} ({\mathfrak X}).$ This functor is compatible with descent, hence we obtain a morphism of stacks over the category of rigid spaces, 

\begin{equation}
\underline D: {\Ccc}_{d, K} \rightarrow  {\mathfrak D}_{d, K} .
\end{equation}

Similarly, if $\mu$ is a miniscule cocharacter with $h_-=0$, $h_+=1$ and with reflex field $E$ as in Remark \ref{mini},
we can define the category $\Ccc_{\mu, K}$
fibered in groupoids over  $E$-rigid spaces by requiring the cokernels $\M/\Phi(\phi^*\M)$ to have a filtration ``of type $\mu$".
In this case, the morphism $\underline D$ sends ${\Ccc}_{\mu, K}$ to  ${\mathfrak D}_{\mu, K} $.

\subsubsection{} It follows from \cite{K2} Theorem (1.2.15) that the functor $\underline D({\rm Sp}(L))$ gives an equivalence ${\Ccc}_{d, K}({\rm Sp}(L))\xrightarrow {\sim}  {\mathfrak D}_{d, K}({\rm Sp}(L))$ for any finite extension $L/\Q_p$.
To briefly explain the construction of the inverse functor we need some notation:
Denote by $\sigma$   the isomorphism   
   $\O\rightarrow \O$  given by applying Frobenius (only) to the coefficients of the power series.
Denote by $x_n$ the point of $\Bbb U$ that corresponds to the irreducible polynomial $E(u^{p^n})$ and let $\hat\O_{{\Bbb U}, x_n}$ the complete local ring of $\Bbb U$ at $x_n$. Notice that the function $\sigma^{-n}(\lambda)\in \O$
has a simple zero at $x_n$. Now consider the composite map
\begin{equation*}
\O\otimes_{K_0}D\xrightarrow {\sigma^{-n}\otimes \Phi^{-n}} 
\O\otimes_{K_0}D\xrightarrow{ } \hat\O_{{\Bbb U}, x_n}\otimes_{K_0}D=
\hat\O_{{\Bbb U}, x_n}\otimes_{K}D_K,
\end{equation*}
where in the first arrow $\Phi^{-n}: D\rightarrow D$ makes sense since $\Phi=\Phi_D$ is bijective. By the above, this induces a map
\begin{equation*}
i_n: \O[\lambda^{-1}]\otimes_{K_0}D\xrightarrow { }
\hat\O_{{\Bbb U}, x_n}[(u-x_n)^{-1}]\otimes_{K}D_K\ .
\end{equation*}
Now suppose we are given 
an object $(D, \Phi, {\rm Fil}^\bullet)$ over $L$. Kisin constructs
a $\Phi$-module 
$(\M, \Phi)=\M(D, \Phi, {\rm Fil}^\bullet)$ 
over $\O_L=\O_{{\rm Sp}(L)\times {\Bbb U}}$ by taking
\begin{equation}\label{intersection}
\M=\bigcap_{n\geq 0} i^{-1}_n( {\rm Fil}^{1}\otimes_{K}(u-x_n)^{-1}\hat\O_{{\Bbb U}, x_n}+D_K\otimes_{K}\hat\O_{{\Bbb U}, x_n})
\end{equation}
and setting $\Phi: \phi^*\M\rightarrow \M$ to be the
restriction of 
$$
1\otimes \Phi_D: \phi^*(\O[\lambda^{-1}]\otimes_{K_0}D)\rightarrow  
\O[\lambda^{-1}]\otimes_{K_0}D\ .
$$
Note that by definition, we have
\begin{equation}
\O\otimes_{K_0}D\subset \M\subset \lambda^{-1}\O\otimes_{K_0}D\ .
\end{equation}
Observe that by its construction, $\M$ is a closed $\O$-submodule of $\lambda^{-1}\O\otimes_{K_0}D$
and so by \cite{K2} Lemma 1.1.4, $\M$ is finite free over $\O$. 

\subsubsection{}\label{5b4} This construction extends to the case that $L$ is a complete rank-$1$
valued field: Let $R^\circ$ 
be a $p$-adic valuation ring of rank $1$ with $L=R^\circ[1/p]$.  Then $\O_L=\O_{[0,1)}$ 
is the ring of rigid functions on the open 
unit disk over   $L\otimes_{\Q_p}K_0$. 
  For simplicity, set  $L_{K_0}=L\otimes_{\Q_p}K_0$,
 $L_K=L\otimes_{\Q_p}K$.
Consider a triple $(D, \Phi, {\rm Fil}^\bullet)$ as above.
We can construct a $\Phi$-module $\M$ over $\O_{[0,1)}$ by (\ref{intersection}) as above  that satisfies
\begin{equation}
\O_{[0,1)}\otimes_{L_{K_0}}D\subset \M\subset \lambda^{-1}\O_{[0,1)}\otimes_{L_{K_0}}D\ .
\end{equation} 
By 
\cite{Gruson}, V,  Rem. $3^\circ$, p. 87, $\O_{[0,1)}$  is a product of Pr\"ufer domains. 
We can see that when an integral power of $r$ is in the set $|L|$ the restriction  $\M_{|[0, r]}$ of $\M$ to the 
closed disk $[0,r]$ is given by a finitely generated torsion free $\O_{[0,r]}$-module
which is free (since $\O_{[0,r]}$ is a product of p.i.ds).
Using \cite{Gruson}, V, Thm. 1, p.\ 83,  we can see that $\M$ is a projective 
finitely generated $\O_{[0,1)}$-module. As such 
it is a direct sum of a free module with a projective module $\L$ of rank $1$
and $\L\simeq {\rm det}(\M)$. We can see that ${\det}(\M)=\lambda^{-a}\O_{[0,1)}$
with $a={\rm dim}_{L_K}({\rm Fil}^1)$; therefore $\L$ and hence $\M$ is finite free
over $\O_{[0,1)}$.

\subsubsection{}  The constructions of the two previous paragraphs are  compatible
in the following sense:  Suppose that the $p$-adic ring $\mathcal A$ is topologically of finite presentation over $\Z_p$. Then $A=\A[1/p]$ is a Tate algebra and we can consider  
the affinoid rigid space ${\rm Sp}(A)$. Recall that ${\rm Sp}(A)$ is the ``generic fiber" ${\rm Spf}(\A)^{\rm rig}$  of  ${\rm Spf}(\A)$ 
in the sense of Raynaud. An object $(\Mfr, \Phi)$ of  
$\widehat \Cc_{d, K}(\A)$ gives an object $(\M, \Phi)$ of 
$\Ccc_{d, K}({\rm Sp}(A))$ by taking $\M$ to be the coherent sheaf
with global sections $\Mfr\otimes_{\A_W[[u]]}\O_A$. This gives a functor
$\widehat \Cc_{d, K}(\A)\rightarrow \Ccc_{d, K}({\rm Sp}(A))$.
The diagram  
$$
\xymatrix{
{} & \ar[d]_{\underline D(\A)}  {\widehat\Cc_{d, K}(\A)}  \ar[r] & \Ccc_{d, K}({\rm Sp}(A))\ar[d]^{\underline D({\rm Sp}(A))}\\
 & {\mathcal D}_{d, K} (A) \ar[r]& {\mathfrak D}_{d, K}({\rm Sp}(A))\ ,}
$$
commutes up to natural equivalence. (Here the lower horizontal arrow is given 
by sending the $A$-module $D$ to the corresponding coherent sheaf over ${\rm Sp}(A)$.) The diagonal   arrow
\begin{equation*}
\Pi(\A): {\widehat\Cc_{d, K}(\A)}\longrightarrow {\mathfrak D}_{d, K}({\rm Sp}(A))
\end{equation*}
(obtained as the composition of the top followed by the right downward arrow)
is by definition, the period functor for $\A$.  It globalizes as follows. 

Suppose that ${\X}$ is an adic formal scheme 
which is locally of finite type over $\Z_p$ (hence $p\,\O_{\mathcal X}$ is an ideal of definition),  and denote by 
${\mathfrak X}={\X}^{\rm rig}$ the corresponding   
rigid space given by its generic fiber, comp. \cite{R-Z}, Prop.\ 5.3.  The construction 
$(\Mfr, \Phi)\mapsto (\M, \Phi)$ above
generalizes to give a functor
$\widehat \Cc_{d, K}({\X})\rightarrow \Ccc_{d, K}({\X}^{\rm rig})$.
Its composition with the functor $\omega({\X^{\rm rig}}):  \Ccc_{d, K}({\X}^{\rm rig})\xrightarrow{ \ }  {\mathfrak D}_{d, K}({\X}^{\rm rig})$
above
gives the  {\it period functor} 
\begin{equation}
\Pi({\X}): \widehat \Cc_{d, K}({\X})\longrightarrow  {\mathfrak D}_{d, K}({\X}^{\rm rig})\ .
\end{equation}
It is localizing in the following sense. Let $\mathcal X=\bigcup\nolimits_i \mathcal U_i$ be an open covering of the formal scheme $\mathcal X$. This induces an admissible open covering of the associated rigid-analytic spaces,
\begin{equation*}
\mathcal X^{\rm rig}=\bigcup\nolimits_i \mathcal U_i^{\rm rig}\ ,
\end{equation*}
comp.\ \cite{R-Z}, Prop.\ 5.3. Then the corresponding diagram of $2$-cartesian rows, with vertical arrows the period morphisms,   is commutative,
$$
\begin{aligned}
 \widehat \Cc_{d, K}&({\X})\rightarrow  &&\prod\nolimits_i\widehat \Cc_{d, K}({\mathcal U_i})&\rightrightarrows &\prod\nolimits_{i, j}  \widehat \Cc_{d, K}({\mathcal U_i\cap \mathcal U_j})\\
   &\downarrow &&\downarrow&\ &\downarrow\\
    {\mathfrak D}_{d, K}&({\X^{\rm rig}})\rightarrow  &&\prod\nolimits_i{\mathfrak D}_{d, K}({\mathcal U_i^{\rm rig}})&\rightrightarrows &\prod\nolimits_{i, j} {\mathfrak D}_{d, K}({\mathcal U_i^{\rm rig}\cap \mathcal U_j^{\rm rig}})\ .
\end{aligned}
$$

\subsubsection{} 
Let us now assume that $R^\circ$ is a complete rank one valuation ring with residue field equal to $\bar\Ff_p$
and set $L=R^\circ[1/p]$ for the corresponding complete rank-$1$ valued field. Recall $\O_L=\O_{[0,1)}$ 
is the ring of rigid functions on the open 
unit disk over $L_{K_0}$. Denote by $\O^{\rm R}_L$ the corresponding
Robba ring
$$
\O^{\rm R}_L:=\lim_{r\rightarrow 1^-} \O_{(r,1)}\ .
$$
This can be identified with the set of Laurent power series $\sum_{n\in \Z}a_n u^n$
with coefficients in $L_{K_0}$ that converge in some open annulus $r<|u|<1$. The ring $\O^{\rm R}_L$
is equipped with a Frobenius endomorphism $\phi: \O^{\rm R}_L\rightarrow \O^{\rm R}_L$
which restricts to $\phi: \O_L\rightarrow \O_L$. 
%Consider also the bounded Robba ring
%$\O^{\rm b,R}_L=\lim_{r\rightarrow 1^-} \O^{\rm b}_{(r,1)}\subset \O^{\rm R}_L$ 
%where $\O^{\rm b}_{(r,1)}$ is the subring of bounded rigid functions. 
Let $\O^{\rm int}_L$ the subring of $\O^{R}_L$ consisting 
of those Laurent power series $\sum_{n\in \Z}a_n u^n$ with $a_n\in R^\circ\otimes_{\Z_p}W$,
for all $n\in \Z$. By \cite{Ke2} Prop. 3.5.5, $\O^{\rm int}_L$ is a henselian local ring with maximal ideal
given by the set of series with $|a_n|<1$, for all $n\in \Z$,
and residue field $\bar{\Bbb F}_p((u))$. Notice that $E(u)$
is a unit in $\O^{\rm int}_L$ and $\phi$ preserves $\O^{\rm int}_L$.

Suppose now that $\N$ is a finite free rank $d$ $\Phi$-module 
over $\O^{\rm R}_L$ with $\Phi:\phi^*\N\xrightarrow{\sim}\N$
an isomorphism. We will say that $(\N,\Phi)$ is {\it purely of slope 
zero} if there is a finite free rank $d$ $\O^{\rm int}_L$-submodule $\N^{\rm int}\subset \N$ such that:

i) $\N^{\rm int}\otimes_{\O^{\rm int}_L}\O^{\rm R}_L= \N$, and

ii) $\Phi_{|\N^{\rm int}}$ induces an isomorphism 
$\phi^*\N^{\rm int}\xrightarrow{\sim}\N^{\rm int}$.
 
When  $L/\Q^{\rm unr}_p$ is finite, this is equivalent to asking that
$(\N, \Phi)$ is purely of slope zero in the sense of Kedlaya \cite{Ke1}.

As in the previous paragraph, we have the period functor
$$
\Pi(R^\circ): \widehat\Cc_{d, K}(R^\circ)\rightarrow {\mathcal D}_{d, K}(L)\simeq {\mathfrak D}_{d,K}({\rm Sp}(L))\ .
$$
 By \ref{5b4}  we can associate to an object $(D, \Phi, {\rm Fil}^\bullet)$ of 
 ${\mathcal D}_{d, K}(L)$ a $\Phi$-module $(\M, \Phi)=\mathcal M(D,\Phi, {\rm Fil}^\bullet)$ over 
 $\O_L$.

\begin{conjecture}
(i) The object $(D, \Phi, {\rm Fil}^\bullet)$ of ${\mathfrak D}_{d, K}({\rm Sp}(L))$ is in the image of the period functor
$\Pi(R^\circ)$,
i.e is of the form $\Pi(R^\circ) (\Mfr, \Phi)$ for some $(\Mfr, \Phi)\in \widehat\Cc_{d, K}(R^\circ)$, if and only if 
$\mathcal M(D,\Phi, {\rm Fil}^\bullet)\otimes_{\O_L}\O^{\rm R}_L$ is purely of slope zero.

%(ii) Any two objects $(\Mfr, \Phi), (\Mfr', \Phi')$ of $\widehat\Cc_{d, K}(R^\circ)$ have isomorphic images under $\Pi$ if and only if they have the same image in 
%$\widehat{\mathcal  R}_d(R)$. 

(ii) There exists an $\Q_p$-analytic subspace (in the sense of Berkovich)
$$
({\rm Res}_{K_0/\Q_p}GL_d\times_{\Q_p} {\rm Gr}_{d,  K})^{\rm an}\subset {\rm Res}_{K_0/\Q_p}GL_d\times_{\Q_p} {\rm Gr}_{d, K}
$$
 invariant under ${\rm Res}_{K_0/\Q_p}GL_d$ such that the fiber over $L$ of the stack quotient 
 $$
 [({\rm Res}_{K_0/\Q_p}GL_d\times_{\Q_p} {\rm Gr}_{d, K})^{\rm an}/_{(\phi, \cdot )} {\rm Res}_{K_0/\Q_p}GL_d ]
 $$
  parametrizes the image points of $\Pi$ over $L$.
\end{conjecture}

There is an obvious variant of this conjecture involving a minuscule cocharacter $\mu$ with $h_{-}=0, h_{+}=1$. 
    
\section{Kisin varieties and Bruhat-Tits buildings}\label{4}

\subsection{}
We return to the set-up and notations of \S \ref{schubert}. Let $\Ff$ be a finite extension of $\Ff_p$. 

\subsubsection{} For simplicity, set $L=\Ff\otimes k$. Recall $G={\rm Res}_{W/\Z_p}GL_d$.
Suppose now that $A\in G(\Ff((u)))=GL_d(L((u)))$ and consider the corresponding
$L((u))$-$\Phi$-module $M_A=(L((u))^d, A\cdot \phi)$ which gives an object of $\RR(\Ff)$.
Choices  $A$, $A'$ that are $\phi$-conjugate, i.e $A'=g^{-1}\cdot A\cdot \phi(g)$  with $g\in GL_d(L((u)))$, give isomorphic modules.
By the above, the fiber product $\{M_A\}\times_{\cal R}{\cal C}_\nu$
is represented over $\Ff$ by a projective subscheme of the affine Grassmannian $\F_G$
of $L[[u]]$-lattices in $L((u))^d$. We denote this subscheme
 by $\Cc_\nu(A)$.
Similarly, we can consider  the fiber product $\{M_A\}\times_{\cal R}{\cal C}^0_\mu$
which is a locally closed subscheme of $\Cc_\nu(A)$; we denote this subscheme by $\Cc^0_\nu(A)$.
This can be thought of as an inseparable analogue of an affine Deligne-Lusztig variety. We call  $\Cc^0_\nu(A)$ the {\it Kisin variety} associated to $(G, A, \nu)$, and   $\Cc_\nu(A)$ the corresponding {\it closed Kisin variety}.

Concretely, for every finite extension $\Ff'$ of $\Ff$, the $\Ff'$-points of the Kisin variety $\Cc^0_\nu(A)$
are given by
$$
\Cc^0_\nu(A)(\Ff')=\{ g\cdot (\Ff'\otimes k) [[u]]^d\ |\  g^{-1}\cdot A\cdot \phi(g)\in G({\Ff'}[[u]])\cdot u^\nu\cdot G({\Ff'}[[u]])\}.
$$  
 The
points $\Cc_{\nu}(A)(\Ff')$ parametrize finite flat commutative group schemes $\Gg$ with $\Ff'$-action  over $\O_K$
which have ``Hodge type $\leq \nu$" and are 
such that the restriction of the Galois representation ${\rm Gal}(\bar K/K)\rightarrow {\rm Aut}_{\Ff'}(\Gg(\bar\O_K))$
to the Galois group $G_\infty$ corresponds to the $\Phi$-module given by $A$.

\subsubsection{} The above construction extends to the set-up of 
a general reductive group $G={\rm Res}_{W/\Z_p}H$ described in \S \ref{genG}.
If  $\nu$ is a dominant coweight 
of $G$ and $A\in LG(\Ff)=G(\Ff((u)))$, we can define as above ${\cal C}^0_{\nu, G}$, ${\cal C}_{\nu, G}$ and the Kisin variety 
${\cal C}^0_{\nu, G}(A)$, and closed Kisin variety ${\cal C}_{\nu, G}(A)$.
However, the relation with Galois representations of ${\rm Gal}(\bar K/K)$ or finite group
schemes is not so clear in this general case.

\subsection{} We now explain how the Bruhat-Tits building can help to get an overview
of a Kisin variety.

\subsubsection{}
For simplicity, we assume that $k=\Ff_p$, $W=\Z_p$ and that $H=G$  is a split Chevalley group over $\Z_p$. In the rest of this section, the symbol $W$ is free again, and will be reserved for Weyl groups. 
 Let $T$ be a maximal split torus
of $G$. We will identify the cocharacter
groups $X_*=X_*(T_{\Q_p})=X_*(T_{\Ff_p})=X_*(T_{\Ff_p((u))})$. Suppose that $C$ is a choice of a positive closed Weyl chamber
in the vector space
$$
V=X_*(T)\otimes_\Z{\Bbb R}. 
$$

Let ${\cal B}={\cal B}(\Ff((u))$ be the Bruhat-Tits building of $G$ over $\Ff((u))$.
This is a metric space with equivariant distance function $d:{\cal B}\times {\cal B}\rightarrow {\Bbb R}$.
We have  the  ``refined" Weyl distance function
$
\delta: {\cal B}\times {\cal B}\rightarrow C
$
which is defined as follows, cf. \cite{KLM}, section 5.1:
Let $x$, $y\in {\cal B}$ and suppose that ${\cal A}$ is an apartment that contains 
both $x$ and $y$. Let $\delta_{\cal A}(x, y)$ be the unique representative in $C$
of the vector $y-x\in V$ and set
$
\delta(x,y)=\delta_{\cal A}(x,y)
$.
(This is independent of the choice of apartment $\cal A$.)
The function $\delta$ is translation $G$-equivariant and 
satisfies the triangle inequality, cf.\cite{KLM}, Remark 3.33, (ii), 
\begin{equation}
\delta(x,z)\leq \delta(x,y)+\delta(y,z)
\end{equation}
for the  order that extends the usual order on dominant coweights. Also,
\begin{equation}
\delta(x,y)=\delta(y,x)^*\ ,
\end{equation}
where $v\mapsto v^*=w_0(-v)$ is the usual involution of $C$ defined by the longest element $w_0$ of the finite Weyl group $W$. 
\subsubsection{} Consider now the homomorphism $\phi: \Ff((u))\rightarrow \Ff((u))$, given $\phi(a)=a$
if $a\in \Ff$, $\phi(u)=u^p$. We will show that it induces a map $\phi: {\cal B}\rightarrow {\cal B}$
with the following properties:
\begin{enumerate}
\item the image of any apartment under $\phi$ is  an apartment, 
\item we have $\phi(g)\cdot \phi(x)=\phi(g\cdot x)$ for any $g\in G(\Ff((u)))$, $x\in {\cal B}$.
\item For $x$, $y\in {\cal B}$, we have
\begin{equation*}\label{homo}
d(\phi(x),\phi(y))=p\cdot d(x,y), \quad \delta(\phi(x),\phi(y))=p\cdot \delta(x,y).
\end{equation*}
\item The map $\phi: {\cal B}\rightarrow {\cal B}$ takes maps geodesics to geodesics;
i.e.,  if $[x,y]\subset {\cal B}$ is the geodesic in ${\cal B}$ joining $x$ and $y$, then 
the image $\phi([x, y])$ is the geodesic $[\phi(x), \phi(y)]$ joining $\phi(x)$ and $\phi(y)$.
\item The map $\phi$ has a unique fixed point, i.e.,  there is a unique $y_0\in \B$
such that $\phi(y_0)=y_0$. The point $y_0$ is a special vertex in $\B$.
\end{enumerate}
 
 Indeed, consider the vertex $y_0$ of $\B$ which is fixed under the subgroup $G(\Ff[[u]])$.
Let $\A_0$ be the apartment in $\Bb$ that corresponds to a constant maximal torus
$T=T_0\otimes_{\Ff}\Ff((u))$ with $T_0\subset G_{\Ff}$; then $y_0$ belongs to 
$\A_0$ and this choice of base point allows us to identify the affine space $\A_0$
with $V=X_*(T)\otimes_\Z{\Bbb R}$. Scaling by $p$ on $V$ now gives 
a well-defined map $\phi_0: \A_0\rightarrow \A_0$ such that $\phi_0(y_0)=y_0$
and which satisfies $\phi_0(n\cdot y)=\phi(n)\cdot \phi_0(y)$ for each $y\in \A_0$
and $n$ in the normalizer $N(T)\subset G(\Ff((u)))$. Now recall that the building $\Bb$
can be described as the quotient of $G(\Ff((u))\times \A_0$ via the equivalence relation
$(g, x)\sim (x', g')$ if there is $n\in N(T)$ such that $x'=n\cdot x$, $g'=ngn^{-1}$.
Using the above we immediately see that $\phi(x, g)=(\phi_0(x), \phi(g))$ respects the equivalence 
relation and gives $\phi: \Bb\rightarrow \B$; since each apartment of $\B$ is of the form $g\cdot \A_0$
we see that the image of an apartment by $\phi$ is also an apartment. The desired properties now follow
easily by using the above and the fact that  any two points $x$, $y\in \B$ are contained in
some apartment and that the geodesic $[x,y]$ is the straight line segment connecting 
$x$ and $y$ in that apartment. Note that the equality $d(\phi(x),\phi(y))=p\cdot d(x,y)$ 
implies that there is at most one fixed point which then has to be the vertex $y_0$ given above.

\smallskip

Note that, by construction, the  group  
$L^+G(\Ff)=G(\Ff[[u]])$ is the stabilizer of $y_0$
in $G(\Ff((u)))$.
The map
\begin{equation}
\iota:\ \F_G(\Ff)=G(\Ff((u)))/G(\Ff[[u]])\hookrightarrow {\cal B}, \quad g\cdot G(\Ff[[u]])\mapsto g\cdot y_0
\end{equation}
allows us to identify the $\Ff$-valued points of the affine Grassmannian with 
a subset of the vertices in the building.

\subsubsection{} 
Suppose now that $A$ is in $G(\Ff((u)))$
and gives an object in $\RR_G(\Ff)$. Then we have a
map $\Phi=A\cdot \phi: {\cal B}\rightarrow {\cal B}$
which also satisfies 
\begin{equation}\label{homo2}
d(\Phi_A(x), \Phi_A(y))=p\cdot d(x,y), \quad \delta(\Phi_A(x),\Phi_A(y))=p\cdot \delta(x,y).
\end{equation}
Then the $\Ff$-valued points of $\Cc_\nu(A)\subset \F_G(\Ff)$ correspond to
the following subset of vertices of the building,
$$
\Cc_\nu(A)=\{x\  \hbox{\rm vertex in $\B\mid x\in {\rm Im}\ \iota$,  $0\leq \delta(x, \Phi_A(x))\leq \nu$ }\}.
$$
If $N/\Ff((u))$ is a finite separable extension, we have an isometric embedding
$$
{\cal B}\hookrightarrow {\cal B}(N).
$$
We will use this to identify 
${\cal B}$ with a subspace of ${\cal B}(N)$.
The map $\Phi_A$ extends to a map ${\cal B}(N)\rightarrow {\cal B}(N)$. 

\begin{prop}\label{fixed} There is a  finite separable extension $M/\Ff((u))$ 
such that the above map $\Phi_A: {\cal B}(M)\rightarrow{\cal B}(M)$ has a  fixed point.
This fixed point is unique in $\cup_{N/\Ff((u))}{\cal B}(N)$.
\end{prop}

\begin{proof}
The uniqueness follows easily from (\ref{homo2}).
For simplicity, we will write $\Phi$ instead of $\Phi_A$.
Consider the ``Lang isogeny"
$$
G_{\Ff((u))}\rightarrow G_{\Ff((u))}\ ;\quad g\mapsto g^{-1}\phi(g).
$$
This is a finite \'etale surjective morphism;
therefore, if $A$ is in $G(\Ff((u)))$, then there
is $g\in G(M)$ for some finite Galois extension  $M/\Ff((u))$
such that $A=g^{-1}\phi(g)$. Consider $x_0=g^{-1}\cdot y_0$ which is a special vertex in ${\cal B}(M)$.  We have 
$$
\Phi(x_0)=A\cdot \phi(g^{-1}\cdot y_0)=g^{-1}\cdot \phi(g)\cdot 
\phi(g)^{-1}\cdot y_0=g^{-1}\cdot y_0=x_0\ 
$$
and so $x_0$ is a fixed point.
\end{proof}
\smallskip

In fact, if $\sigma$ 
is an element of ${\rm Gal}(M/\Ff((u)))$, since $\sigma\cdot \phi=\phi\cdot\sigma$, we can see that
$\sigma(A)=\sigma(g^{-1}\phi(g))=\sigma(g)^{-1}\phi(\sigma(g))=A=g^{-1}\phi(g)$;
therefore $\sigma(g)g^{-1}\in G(\Ff)$.
Since $g_0\cdot y_0=y_0$ for $g_0\in G(\Ff)$, 
the point $x_0=g^{-1}\cdot y_0$   depends only on $A$
and is ${\rm Gal}(M/\Ff((u)))$-fixed.
Therefore, if $M/\Ff((u))$ is tamely ramified,
 which implies ${\cal B}(M)^{{\rm Gal}(M/\Ff((u)))}={\cal B}$, we conclude that $x_0$
 belongs to ${\cal B}$.

\smallskip
 
\subsubsection{} We continue to write $\Phi=\Phi_A$.
If $x$ is in ${\cal B}$, we can apply the triangle inequality above to $x$, $\Phi(x)$ and $x_0=\Phi(x_0)$, in two different ways.
We obtain:
\begin{equation*}
\delta(\Phi(x), x) \leq \delta(\Phi(x), x_0)+\delta(x_0, x)=p\cdot \delta(x, x_0)+\delta(x, x_0)^*,
\end{equation*}
\begin{equation*}
\delta( \Phi(x), x_0)=p\cdot \delta(x, x_0)\leq \delta(\Phi(x), x)+\delta(x , x_0).
\end{equation*}

Combining these we get 
\begin{equation}\label{ineq}
(p-1)\cdot \delta(x, x_0)\leq \delta(\Phi(x), x)\leq p\cdot \delta(x, x_0)+\delta(x, x_0)^*.
\end{equation}
This implies that if $h\in G(\Ff((u)))$ is such that
\begin{equation}
p\cdot \delta(h\cdot y_0, x_0)+\delta(h\cdot y_0, x_0)^*\leq \nu,
\end{equation}
 then the corresponding point $h\cdot G(\Ff[[u]])$ in $\F_G(\Ff)$ belongs to $\Cc_{\nu, G}(A)$, which is then non-empty and is contained in the ball of radius $\nu/(p-1)$ around $x_0$.

\subsubsection{} Suppose that $A'=h^{-1}\cdot A\cdot \phi(h)$ with $h\in G(\Ff((u)))$.
Then $A'=(gh)^{-1}  \phi(gh)$ 
and the corresponding $\Phi_{A'}$-fixed vertex is $x'_0=(gh)^{-1}\cdot y_0=h^{-1}\cdot x_0$.
We conclude that the orbit $G(\Ff((u)))\cdot x_0$  
 only depends on the $\phi$-conjugacy class
 of $A$ in $G(\Ff((u)))$. By the above, if $M/\Ff((u))$ is tamely ramified,
   $x_0$
 belongs to ${\cal B}$.
 \begin{comment}
We should be able to show (???) that
 for all $x\in {\cal B}$ we have: 
 $$
 \delta(x, x_0)=\delta(x, \tilde x_0)+\delta(\tilde x_0, x_0).
 $$
 This should help in exploiting the inequalities (\ref{ineq}).
 \end{comment}

\subsection{}
We continue to assume that $k=\Ff_p$  and now take $G=H=GL_d$. Take $T$
the standard maximal torus of $GL_d$.  
Then  the finite 
Weyl group is the symmetric group $S_d$, $X_*(T)_\R=X_*(T)\otimes_\Z\Bbb R$, and the standard
choice of a positive closed Weyl chamber is
$$
C=\{(v_1,\ldots , v_d)\in {\Bbb R}^d\ |\ v_1\geq v_2\geq\cdots \geq v_d\}.
$$
The partial  order on $C$ is given by: $(v_1,\ldots, v_d)\leq (v'_1,\ldots ,v'_d)$ iff   
$$
\sum_{i=1}^{r}v_i\leq \sum_{i=1}^{r}v'_i, \hbox{\rm \ for $r=1,\ldots , d-1$, and\ } v_1+\cdots +v_d=v'_1+\cdots +v'_d.
$$
In this case, we will explain the construction of the fixed point in a slightly different way. Start with $M_A=(k((u))^d, A\cdot\phi)$ and
set
$$
U=(k((u))^{\rm sep}\otimes_{k((u))}M_A)^{\phi\otimes\Phi_A={\rm Id}}
\subset  k((u))^{\rm sep}\otimes_{k((u))}M_A 
$$ 
for the  $k$-vector space of the corresponding 
${\rm Gal}(k((u))^{\rm sep}/k((u)))$-representation $\rho$.
(Here $\phi: k((u))^{\rm sep}\rightarrow k((u))^{\rm sep}$ denotes again the Frobenius of the separable closure.) In fact, one can see from the construction of $\rho$
that there is a finite separable extension $L/k((u))$
such that 
$$
U=(L\otimes_{k((u))}M_A)^{\phi\otimes\Phi ={\rm Id}}
\subset L\otimes_{k((u))}M_A=L^d
$$
as ${\rm Gal}(k((u))^{\rm sep}/k((u)))$-modules.
(Note that $L^d$ also supports a $\Phi$-module
structure for the extension $\phi_{|L}$ of $\phi$ to $L$;
this is 
given by $A\cdot \phi_{|L}$
)
Now set ${\frak M}_0$ for the $\O_L$-submodule in 
$L^d$ generated by the elements in $U$. Then ${\frak M}_0/ u_L{\frak M}_0
\simeq U$ and so ${\frak M}_0$ is an $\O_L$-lattice in $L^d$.
Since $\phi\otimes\Phi=A\cdot \phi_{|L}$ acts as identity on $U$, we can see
that 
$$
(A\cdot \phi_{|L})^*({\frak M}_0)={\frak M}_0.
$$ 
The lattice ${\frak M}_0$ gives a point $x_0$ of ${\cal B}(L)$
which is fixed under the map $\Phi$, i.e $\Phi(x_0)=x_0$.
\bigskip

\subsection{} 
In this paragraph, we will explain the picture in the building for $\Ff=\Ff_p$ and $G=H=GL_2$.
Our main objective is the following. Given a dominant coweight $\nu=(a, b)$ with $a\geq b\geq 0$
and a matrix $A\in GL_2(\Ff((u)))$,  describe the set of vertices in the building $\B$ 
which correspond to $\Ff$-valued points in $\Cc_\nu(A)$, i.e, to lattices $\Mfr\subset \Ff((u))^2$
for which $\Phi_A(\phi^*(\Mfr))\subset \Mfr$ and such that  $\Mfr/\Phi_A(\phi^*(\Mfr))=\Mfr/\langle A\cdot \phi(\Mfr)\rangle$ has elementary divisors
$(a', b')$, $a'\geq b'\geq 0$ which are smaller than $\nu=(a, b)$, i.e $a'\leq a$, $a'+b'=a+b$.
The corresponding set in the building  is the set  of vertices $x$ 
 such that $0\leq\delta(x, \Phi_A(x))\leq \nu$.
To simplify our discussion, we will consider the projection $\B\rightarrow \T$
where $\T$ is the tree of homothety classes of lattices in $\Ff((u))^2$
(i.e the building for $PGL_2(\Ff((u)))$).
Note that the Weyl chamber distance $\delta$ on the tree $\T$ coincides (up to sign) 
with the usual distance $d: \T\times \T\rightarrow {\Bbb R}_{\geq 0}$.

We consider the sets
${\rm Vert}(\T)_{\nu, A}$ of vertices $x$ in the tree $\T$ for which
$$
 d(x, \Phi_A(x))\leq r=|a-b|.
 $$
Let $x_0$ be the fixed point of $\Phi=\Phi_A$ on $\T(M)$ and $\ti x_0$ its projection to $\T$.
Note that the inequalities (\ref{ineq}) imply that 
$$
B(x_0, \frac{r }{p+1})\cap {\rm Vert}(\T)\subset  {\rm Vert}(\T)_{\nu, A} \subset B(x_0, \frac{r }{p-1})\cap {\rm Vert}(\T)
$$
where  $B(x_0, R)$ is the ``ball" 
\begin{equation}
d(x, x_0)\leq R
\end{equation}
of radius $R$ centered at the point $x_0$.

To refine this, we will consider several possible cases:

A) $\ti x_0$ is not a vertex in $\T$. Then $\ti x_0$ lies on a segment $[\eta, \eta']$
with $\e$, $\e'$ the closest vertices to $\ti x_0$. Now consider the images $\Phi(\e)$, 
$\Phi(\e')$. Since the geodesic $[x_0, \Phi(\e)]$ passes through the projection $\ti x_0$,
it also has to pass through either $\e$ or $\e'$ (but not both). Similarly for $[x_0, \Phi(\e')]$.
There are several subcases:
\smallskip

1) $\e\in [x_0, \Phi(\e')]$, $\e'\in [x_0, \Phi(\e)]$. Apply $\Phi$ to conclude that $\Phi(\e)$ lies in the geodesic
from $x_0$ to $\Phi^2(\eta')$ and $\Phi(\e')$ in the geodesic from $x_0$ to $\Phi^2(\eta)$.
Note that if $\Phi(\ti x_0)\neq \ti x_0$ and is, for example, between $\ti x_0$ and $\e'$, then 
$\Phi([x_0, \e'])=[x_0, \Phi(\e')]$ would pass first through $\ti x_0$, then through $\Phi(\ti x_0)$,
and then through $\eta$. This contradicts the fact that this is a geodesic. A similar contradiction is obtained
if  $\Phi(\ti x_0)$   is between $\ti x_0$ and $\e$. We conclude
that $\Phi(\ti x_0)=\ti x_0$ and hence $x_0=\ti x_0$.

\begin{figure}[h]
\includegraphics[width=8cm]{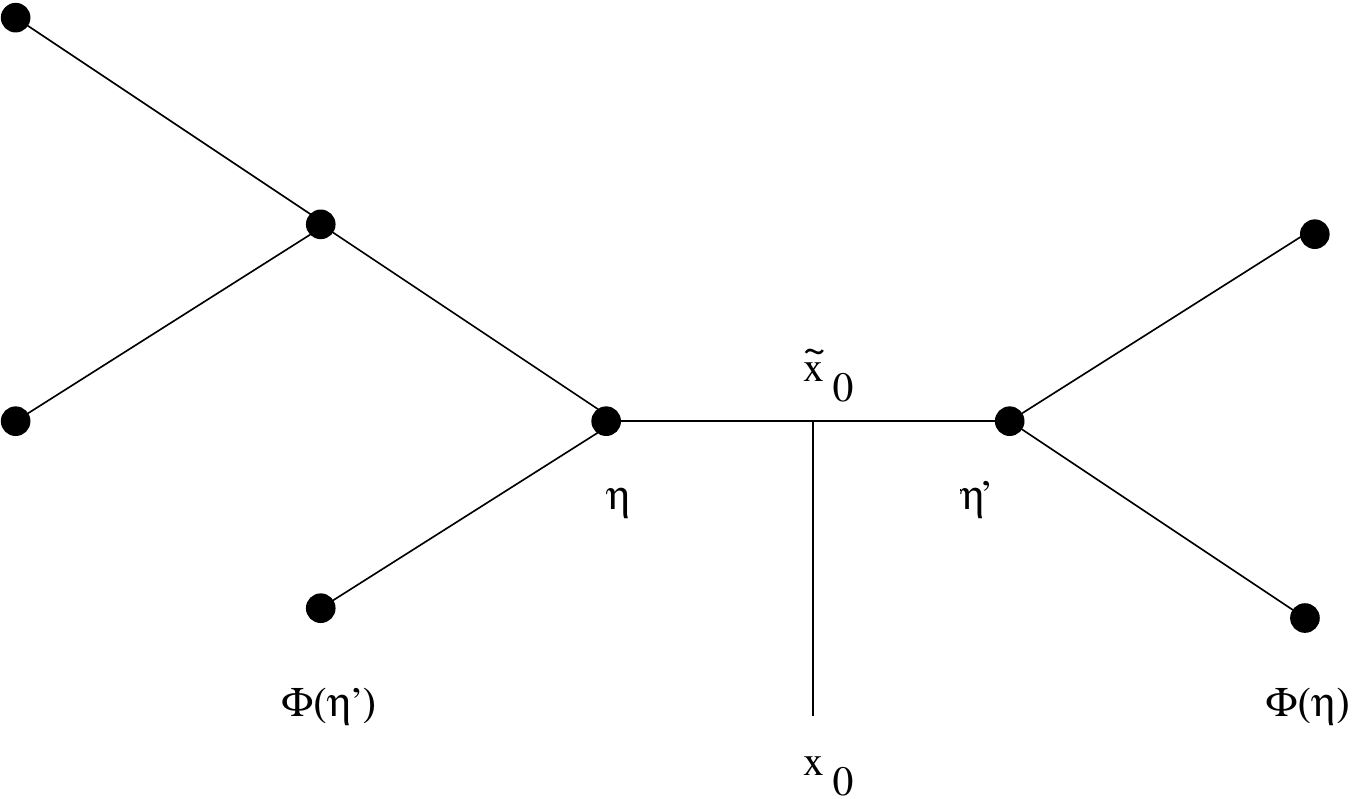}
\caption{The case A1}
\end{figure}

By similar arguments, we deduce that 
the limit 
$$
\ds{\lim_{n\rightarrow \infty}[\Phi^{2n}(\e), \Phi^{2n}(\e')]=\lim_{n\rightarrow \infty}[\Phi^{2n+1}(\e'), \Phi^{2n+1}(\e)]}
$$
gives an apartment which is preserved (but flipped) by $\Phi$.
Indeed, $\Phi$ takes the half-apartment $\ds{\lim_{n\rightarrow \infty}[x_0, \Phi^{2n}(\e)]}$ to 
$\ds{\lim_{n\rightarrow \infty}[x_0, \Phi^{2n+1}(\e)]}$. 

Note that in this case, there are no half-apartments in the
tree $\T$ that are preserved by $\Phi$. Indeed, consider a vertex $y$ in 
such a half-apartment and connect this to $x_0$;
the geodesic has to pass through either $\e$ or $\e'$;
in either case, since the geodesic from $x_0$ to $\Phi(y)$
has to pass through the opposite point $\e'$, resp. $\e$,
we obtain a contradiction.
Recall that the set of half apartments in the tree can be naturally identified
with the set of one-dimensional subspaces of the corresponding vector space
$\Ff((u))^2$. Hence, we see that in  case (A1) 
 the $\Phi$-module given by the matrix $A$ is simple.

Now suppose $x$ is a vertex in $\T$. The geodesic 
$[x_0, x]$ has to pass through either $\e$ or $\e'$. Suppose that $\e'\in [x_0, x]$ (the other case is similar) and consider
$\Phi([x_0, x])=[\Phi(x_0), \Phi(x)]=[x_0, \Phi(x)]$. This contains $\Phi(\e')$ and therefore has to pass
through $\e$ (since $\e\in [x_0, \Phi(\e')]$). Therefore, the geodesic $[x, \Phi(x)]$ passes through both $\e$ and $\e'$
(and also $ x_0$) and we have 
\begin{eqnarray*}
d(x, \Phi(x))&=&d(x,  x_0)+d( x_0, \Phi(x))\ \ \ \ \ \\
&=&(p+1)d(x, x_0) .
\end{eqnarray*}
Hence, in this case, $d(x, \Phi_A(x))\leq r$ amounts to $d(x, x_0)\leq r/(p+1)$
and we have
\begin{equation*}
{\rm Vert}(\T)_{\nu, A}=B(x_0, \frac{r}{p+1})\cap {\rm Vert}(\T).
\end{equation*}
 
\smallskip

2) $\e\in [x_0, \Phi(\e)]$, $\e'\in [x_0, \Phi(\e')]$. Then, we can see that the limits $\lim_{\rightarrow n}[x_0, \Phi^n(\e)]$,
 $\lim_{\rightarrow n}[x_0, \Phi^n(\e')]$ give two half-apartments that are both 
preserved by $\Phi$. As above, this implies 
 that the $\Phi$-module given by $A$ contains two
 $1$-dimensional $\Phi$-submodules; we can easily see that these are distinct. Hence the $\Phi$-module given by $A$ is decomposable.

\begin{figure}[h]
\includegraphics[width=8cm]{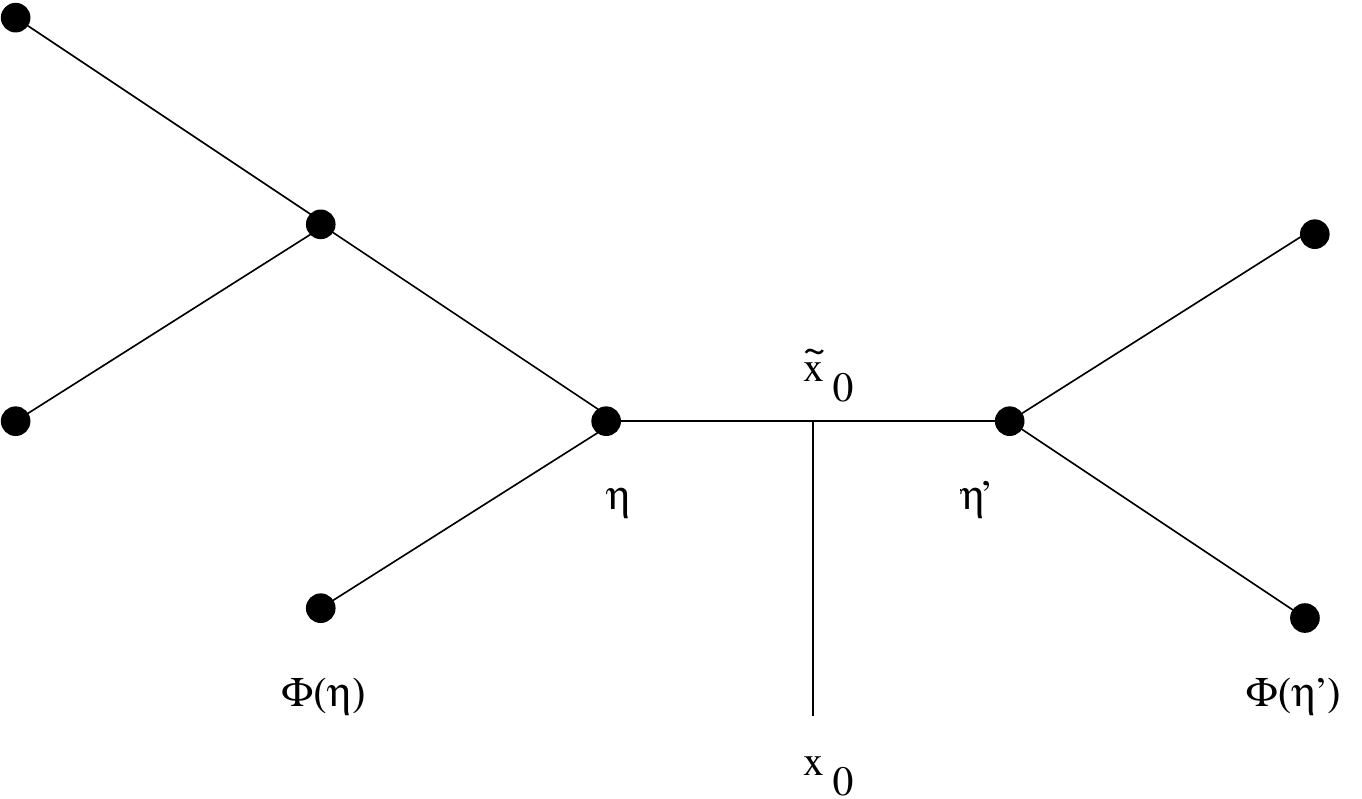}
\caption{The case A2}
\end{figure}

 Now suppose $x$ is a vertex in $\T$. The geodesic 
$[x_0, x]$ has to pass through either $\e$ or $\e'$. Suppose that $\e\in [x_0, x]$
and in fact suppose that $a\geq 0$ is the largest integer for which $\Phi^a(\e)$
is contained in $[x_0, x]$. Consider $\Phi([x_0, x])=[x_0, \Phi(x)]$ which has to contain $\Phi^{a+1}(\e)$.
Therefore, the geodesic $[x, \Phi(x)]$ has to pass through $\Phi^{a+1}(\e)$.
We obtain
\begin{eqnarray*}
d(x, \Phi(x))&=&d(x,  \Phi^{a+1}(\e))+d(  \Phi^{a+1}(\e), \Phi(x))\ \\
&=&d(x,  \Phi^{a+1}(\e))+p\cdot d(  \Phi^a(\e),  x).
\end{eqnarray*}
If $x'$ is the projection of $x$ to the half-apartment $\lim_{\rightarrow n}[x_0, \Phi^n(\e)]$,
then we can rewrite this distance as
\begin{eqnarray*}
d(x, \Phi(x))&=&d(\Phi(x),  \Phi(x'))+d(  \Phi(x'), \Phi^{a+1}(\e))+ \\
&&+d(\Phi^{a+1}(\e), \Phi^a(\e))-d(x', \Phi^a(\e))+d(x,x')\\
&=&(p+1)d(x, x')+(p-1)d(x', \Phi^a(\e))+p^a d(\Phi(\e), \e).
\end{eqnarray*}
There is a similar expression if $\e'$ is in $[x_0, x]$.
Hence, $d(x, \Phi_A(x))\leq r$ can be described as the union of 
two ``thinning tubes" around the two half-apartments that are preserved 
by $\Phi$. Note that, in the above, when $d(x, \Phi(x))$ is bounded, the possible values
of $a$ are bounded too.

3) $\e\in [x_0, \Phi(\e)]$, $\e\in [x_0, \Phi(\e')]$ (the case $\e'\in [x_0, \Phi(\e')]$, $\e'\in [x_0, \Phi(\e)]$ is symmetric). 
Then $\lim_{\rightarrow n}[x_0, \Phi^n(\e)]$
  gives a half-apartment which is preserved 
 by $\Phi$. As above, this implies 
 that the $\Phi$-module given by $A$ contains a
 $1$-dimensional $\Phi$-submodule
 and, therefore, it is not simple. In this case, we can see
 that this is the unique half-apartment preserved by $\Phi$.
 Indeed,  consider a vertex $y$ in 
such a half-apartment $\A'$ and connect this with a geodesic to $x_0$;
the geodesic has to pass through either $\e$ or $\e'$ and we can easily rule out $\e'$.
Now the geodesic $[x_0, \Phi(y)]$ has to pass through $\Phi(\e)$.
Since $\Phi(y)$ is also in $\A'$, we can conclude that $\A'$ 
also contain $\Phi(\e)$. Inductively, $\A'$ contains $\Phi^n(\e)$
for all $n$. We conclude that the $\Phi$-module given by $A$ is not simple and not  decomposable.

\begin{figure}[h]
\includegraphics[width=8cm]{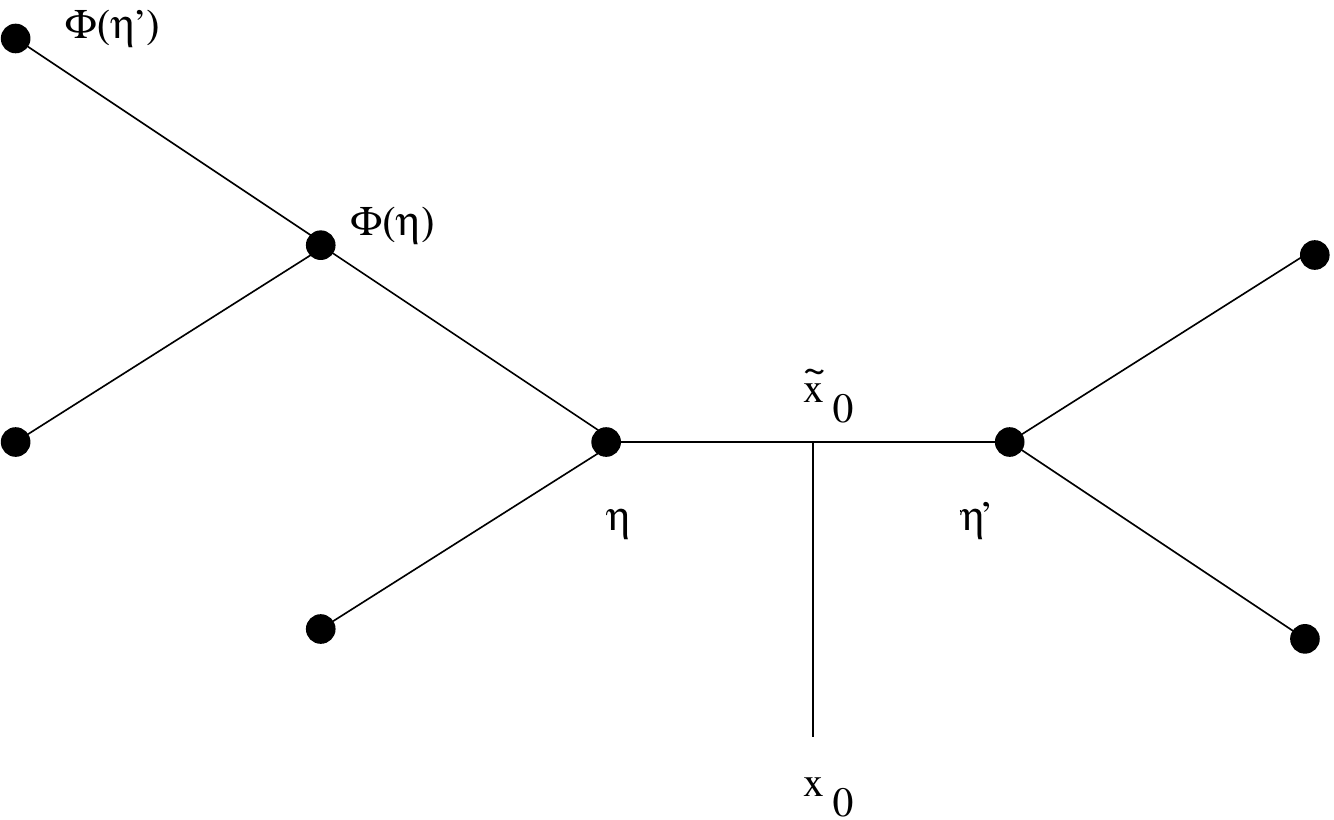}
\caption{The case A3}
\end{figure}

 Now suppose $x$ is a vertex in $\T$. The geodesic 
$[x_0, x]$ has to pass through either $\e$ or $\e'$. 

$\bullet$ Suppose first that the geodesic $[x_0, x]$ passes through $\e$. In fact, suppose that $a\geq 0$ is the largest integer for which $\Phi^a(\e)$
is contained in $[x_0, x]$. Then as before, we obtain
\begin{eqnarray*} 
d(x, \Phi(x))&=&d(x,  \Phi^{a+1}(\e))+pd(  \Phi^a(\e),  x)\\
&=&(p+1)d(x, x')+(p-1)d(x', \Phi^a(\e))+p^a d(\Phi(\e), \e)
\end{eqnarray*}
with $x'$ the projection of $x$ to the apartment $\lim_{\rightarrow n}[x_0, \Phi^n(\e)]$. Hence again the set $d(x, \Phi_A(x))\leq r$ can be described for these vertices as a union of thinning tubes. 

$\bullet$ Now suppose that the geodesic $[x_0, x]$ passes through $\e'$. Then an argument as in case (A1) gives
\begin{equation*} 
d(x, \Phi(x))=(p+1)d(x, x_0)-2d(\ti x_0, x_0).
\end{equation*}
Hence for this kind of vertices this set $d(x, \Phi_A(x))\leq r$ is a ball around $x_0$. 

\smallskip

B) Suppose now that $\ti x_0=\e$ is a vertex of $\T$.
There are two subcases:

1) $\ti x_0=\e$ is not fixed by $\Phi$.  
Then, $\lim_{\rightarrow n}[x_0, \Phi^n(\e)]$
  gives again a half-apartment that is preserved 
 by $\Phi$. We can, in fact, see as before, that this is the unique such half-apartment.
 Hence, in this case, the $\Phi$-module given by $A$ is not simple and not decomposable. Note that after replacing $\Ff$ by a finite extension, $x_0$ becomes of type B2 below.

\begin{figure}[h]
\includegraphics[width=8cm]{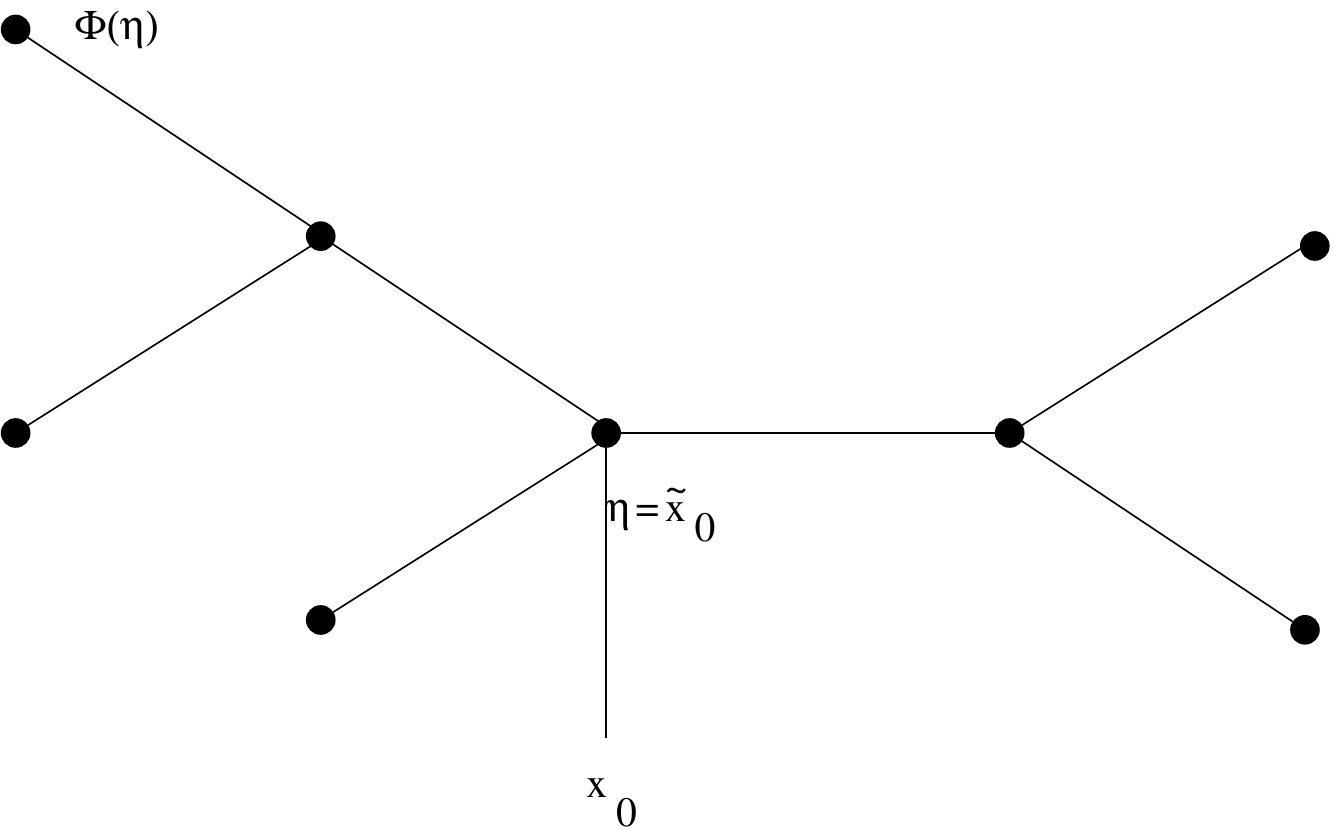}
\caption{The case B1}
\end{figure}

 If $x$ is a vertex of $\T$ then $\e\in [x_0, x]$. If $a\geq 0$ is the largest integer for which $\Phi^a(\e)$
is contained in $[x_0, x]$ we obtain for $d(x, \Phi(x))$ the same formula
as in cases A2 or A3a. Hence in this case $d(x, \Phi_A(x))\leq r$ is a union of thinning tubes. 

2) $\ti x_0=\e$ is fixed by $\Phi$, in other words the fixed point $x_0$ is a vertex of $\T$.
This corresponds to the homothety class of a lattice $\Mfr_0$; we have $\Phi(\phi^*\Mfr_0)=u^s\Mfr_0$ for some
$s$. Denote by $\{\e_i\}_{i=0,\ldots , p}$ its neighborhood vertices. The link of $\eta$
is identified with the projective space of lines in $\Mfr_0/u\Mfr_0$ and the action of $\Phi$ on the link 
then corresponds to the action on the projective space given the linear action
of  $u^{-s}\cdot \Phi$ on   $\Mfr_0/u\Mfr_0$. Now observe that the geodesic $[x_0, \Phi(\eta_i)]$
passes through $\eta_i$ if and only if the action of $\Phi$ on the link leaves the point of the link that is given by $\eta_i$ fixed. Depending on whether the number of fixed points in the link is $\geq 2$, resp.\ $1$, resp.\ $0$, the $\Phi$-module is decomposable, resp. not simple and not decomposable, resp. simple. Note that in the last case it is not absolutely simple, it disappears if $\Ff$ is replaced by a finite extension.

Now let $x$ be a vertex of $\T$ and consider the geodesic  $[x_0, x]$ which has to pass through one of the vertices
$\eta_j$. We distinguish cases according as $\eta_j$ gives a fixed point of the link, or not.

\begin{figure}[h]
\includegraphics[width=8cm]{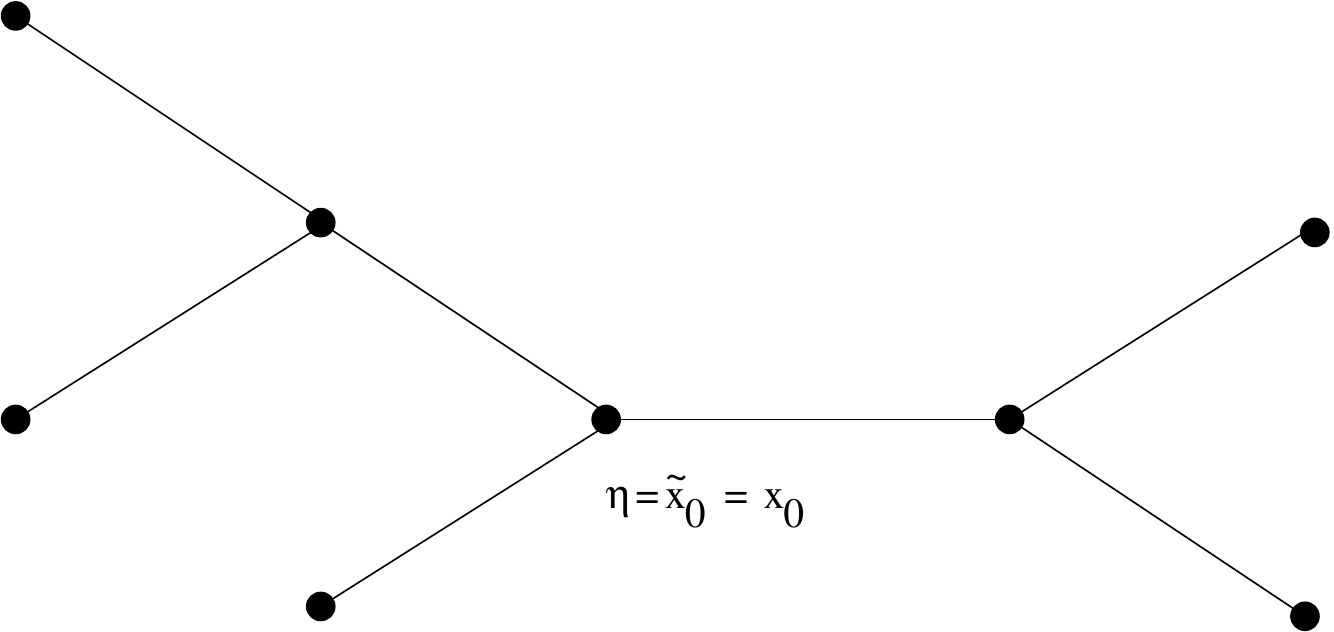}
\caption{The case B2}
\end{figure}

$\bullet$ Suppose first that  $[x_0, x]$ passes through a vertex $\eta_j$ with the corresponding point of the link fixed by $\Phi$.
Then the argument of case A2 applies to obtain $d(x, \Phi(x))$. (It involves the  largest integer $a\geq 0$ such that  
$\Phi^a(\eta_j)$ is in $[x_0, x]$.) For this kind of vertices $d(x, \Phi_A(x))\leq r$ is a union of thinning tubes. 

$\bullet$ Now suppose that $[x_0, x]$ passes through a vertex $\eta_j$ with the corresponding point of the link not fixed by $\Phi$.
Then the argument of case A1 applies to give
\begin{equation*}
d(x, \Phi(x))=(p+1)d(x, x_0).
\end{equation*}
For this kind of vertices  $d(x, \Phi_A(x))\leq r$ is a ball around $x_0$.


\bigskip

 
\bigskip

\begin{thebibliography}{XYZU}


 \bibitem [BL] {BL} A.\ Beauville,  Y.\ Laszlo: Un lemme de descente,   C. R. Acad. Sci. Paris S\'er. I Math. {\bf 320} (1995), no. 3, 335--340.
 
 \bibitem [BG] {BG} S.\ Bosch, U.\ G\"ortz: Coherent modules and their descent on relative rigid spaces,  J. Reine Angew. Math. {\bf 495} (1998), 119--134.

 \bibitem [BTI] {BTI} F.\  Bruhat, J.\ Tits: Groupes r\'eductifs sur un corps local,  Inst. Hautes Etudes Sci. Publ. Math. {\bf 41} (1972), 5--251. 

 \bibitem [BTII] {BTII} F.\  Bruhat, J.\ Tits: Groupes r\'eductifs sur un corps local. II. Sch\'emas en groupes. Existence d'une donn\'ee radicielle valu\'ee,  Inst. Hautes Etudes Sci. Publ. Math. {\bf 60} (1984), 197--376.
 
 \bibitem [Br] {Br} C.\ Breuil: Sch\'emas en groupes et corps des normes, unpublished (1998). 
 
 \bibitem[Ca]{Ca} X.\ Caruso: Sur la classification de quelques $\phi$-modules simples, Preprint 2008. arXiv:0807.1719  
 
 \bibitem[Dr] {Dri}  V.\ Drinfeld: Infinite-dimensional vector bundles in algebraic geometry: an introduction, in: The unity of mathematics, 263--304, Progr. Math., 244, Birkh{\"a}user Boston, Boston, MA, 2006. 
 
  \bibitem[E-K] {emerton-kisin} M.\ Emerton,  M.\ Kisin: Unit $L$-functions and a conjecture of Katz, Ann. of Math. (2) {\bf 153} (2001), no. 2, 329--354. 
  
   \bibitem[Fo] {Fontaine} J-M.\ Fontaine: Repr\'esentations $p$-adiques des corps locaux. I,  The Grothendieck Festschrift, Vol. II, 249--309, Progr. Math., 87, Birkh\"auser Boston, Boston, MA, 1990. 
   
 \bibitem[GHKR] {ghkr}  U.\ G\"ortz, T.\ Haines, R.\ Kottwitz, D.\ Reuman: Dimensions of some affine Deligne-Lusztig varieties, 
 Ann. Sci. \'Ecole Norm. Sup. (4) {\bf 39} (2006), no. 3, 467--511. 
 
  
 \bibitem[Gr] {Gruson} L.\ Gruson: Fibr\'es vectoriels sur un polydiscue
 ultram\'etrique. Ann. Sci. \'Ecole Norm. Sup.   {\bf 1} (1968), no. 1, 45-89.
 
 \bibitem[Ha] {Ha}  U.\ Hartl: On a Conjecture of Rapoport and Zink, arXiv:math/0605254 
 
 \bibitem [He] {He} E.\ Hellmann: On the structure of some moduli spaces of finite flat group schemes, Preprint 2008. arXiv:0810.5277  
 
\bibitem[KLM] {KLM} M.\ Kapovich, B.\ Leeb, J.\ Millson: Convex functions on symmetric spaces, side lengths of polygons and stability inequalities for weighted configurations at infinity, math.DG/0311486

\bibitem[Ka] {katz} N.\ Katz, $p$-adic properties of modular schemes and modular forms, in: Modular Functions of One Variable III ( Proc. Internat. Summer School, Univ. Antwerp, Antwerp, 1972), Lecture Notes in Math. 350, Springer-Verlag, New York, 1973, 69--190. 

\bibitem[Ke1] {Ke1} K.\  Kedlaya: Slope filtrations revisited, Doc. Math. {\bf 10} (2005), 447-525; errata, ibid. {\bf 12} (2007), 361-362.

\bibitem[Ke2] {Ke2} K.\  Kedlaya:  Local monodromy of $p$-adic differential equations: an overview, Int. J.  of Number Theory {\bf 1} (2005), 109-154.

\bibitem[Ki1] {K1} M.\ Kisin:  Moduli of finite flat group schemes and modularity, Annals of Math., to appear.

\bibitem[Ki2] {K2} M.\ Kisin: Crystalline representations and $F$-crystals, in: Algebraic geometry and number theory, 459--496, Progr. Math., 253, Birkh\"auser Boston, Boston, MA, 2006. 

\bibitem[Ki3]{K3} M.\ Kisin: Potentially semi-stable deformation rings, J. Amer. Math. Soc. {\bf 21} (2008), 513--546.

\bibitem[La] {La} Y.\ Laszlo: A non-trivial family of bundles fixed by the square of Frobenius, C. R. Acad. Sci. Paris S\'er. I Math. {\bf 333} (2001), no. 7, 651--656. 

\bibitem[Ma] {Ma} B.\ Mazur: An introduction to the deformation theory of Galois representations, in:  Modular forms and Fermat's last theorem (Boston, MA, 1995), 243--311, Springer, New York, 1997.

\bibitem [PR1] {P-R1} G.\ Pappas, M.\ Rapoport: Local models in the ramified case I. The EL-case,  J. Alg. Geom.  {\bf 12} (2003), 107--145.

\bibitem [PR2] {P-R2} G.\ Pappas, M.\ Rapoport: Local models in the ramified case II. Splitting models,
Duke Math. Journal {\bf 127} (2005), 193--250.

\bibitem [PR3] {P-R3} G.\ Pappas, M.\ Rapoport: Twisted loop groups and their affine flag varieties, Adv. Math. {\bf 219} (2008), 118--188. 


\bibitem[Ra] {Ra} R.\ Ramakrishna:
On a variation of Mazur's deformation functor, 
Compositio Math. {\bf 87} (1993), no. 3, 269--286. 

\bibitem [R] {R} M.\ Rapoport: A guide to the reduction modulo $p$ of Shimura varieties.  Ast\'erisque {\bf 298} (2005), 271--318.

\bibitem [RZ] {R-Z} M.\ Rapoport, Th.\ Zink: {\sl Period spaces for $p$--divisible groups},
Ann.\ of Math. Studies {\bf 141}, Princeton University Press (1996).


\bibitem[T] {T} J.\ Tits: Reductive groups over local fields, in:  {\sl Automorphic forms, representations and $L$-functions.} (Proc. Sympos. Pure Math., Oregon State Univ., Corvallis, Ore., 1977), Part 1,  pp. 29--69, Proc. Sympos. Pure Math., XXXIII, Amer. Math. Soc., Providence, R.I. (1979).

 
 
 
 

\end{thebibliography}
\end{document}